\documentclass{amsart}
\usepackage{graphicx}
\usepackage{amsmath, amssymb, amsfonts}
\usepackage{hyperref}
\usepackage[shortalphabetic]{amsrefs}
\usepackage{mathptmx}    

\begin{document}
\title{Contraction of convex hypersurfaces by nonhomogeneous functions of curvature} 

\author[J. A. McCoy]{James A. McCoy}
\address{Priority Research Centre Computer Assisted Research Mathematics and Applications, School of Mathematical and Physical Sciences, University of Newcastle, University Drive, Callaghan, NSW 2308, Australia}
\email{James.McCoy@newcastle.edu.au}
\keywords{curvature flow, parabolic partial differential equation, hypersurface}
\subjclass[2000]{53C44}

\begin{abstract}
A recent article \cite{LL20} considered contraction of convex hypersurfaces by certain nonhomogeneous functions of curvature, showing convergence to points in finite time in certain cases where the speed is a function of a degree-one homogeneous, concave and inverse concave function of the principle curvatures.  In this article we extend the result to various other cases that are analogous to those considered in other earlier work, and we show that in all cases, where sufficient pinching conditions are assumed on the initial hypersurface, then under suitable rescaling the final point is asymptotically round and convergence is exponential in the $C^\infty$-topology.
\end{abstract}

\keywords{curvature flow, parabolic partial differential equation, hypersurface, mixed volume, axial symmetry}

 \thanks{This research was supported by Discovery Grant DP180100431 of the Australian Research Council.  Part of this work was completed while the author was at the Okinawa Institute of Science and Technology as part of the Visiting Mathematics Professors program.}

\maketitle

\section{Introduction} \label{S:intro}
\newtheorem{Fconditions}{Conditions}[section]
\newtheorem{Phiconditions}[Fconditions]{Conditions}
\newtheorem{main}[Fconditions]{Theorem}
The contraction of convex hypersurfaces by their curvature in Euclidean space has been very well studied since Huisken's seminal work on the mean curvature flow \cite{Hu}.  For a review of many fully nonlinear contraction flows in Euclidean space and their action on convex initial hypersurfaces we refer the reader to \cite{AMZ11} and the references contained therein.  A feature all the works here is that the normal flow speeds are homogeneous functions of the principal curvatures.  Homogeneity is used in several places in the argument to prove convergence to round points, often in conjunction with Euler's identity for such functions.

Contraction flows by speeds that are not homogeneous have been much less well-studied.  Alessandroni-Sinestrari considered asymptotic behaviour of a flow by a particular function of mean curvature \cite{AS}; this speed was later investigated in the setting of convex hypersurfaces by Espin \cite{E}.  Smoczyk considered flow by functions of the mean curvature in \cite{Sm97}.  Chou and Wang considered a kind of logarithmic anisotropic Gauss curvature flow for convex hypersurfaces in \cite{CW}.  For expansion flows, Chow and Tsai considered nonhomogeneous Gauss curvature flows in \cite{CT}, proving asymptotic roundness for expansion flows by functions of Gauss curvature with certain properties and some partial results for contraction and bidirectional flows.  They considered expansion by more general nonhomogeneous speeds in \cite{CTasian}.    

Some general results for flowing surfaces by nonhomogeneous speeds were obtained in \cite{GK}.  Very recently the work of Li and Lv considered contraction of convex hypersurfaces by nonhomogeneous speeds \cite{LL20}.  Specifically their speeds are a function $\Phi$ of a degree-one homogeneous function of the principal curvatures $F$ where $F$ and $\Phi$ each satisfy some natural properties, in particular, $F$ is concave and inverse concave.  The authors showed that strictly convex initial hypersurfaces shrink in finite time to points.  (In their paper they also obtained other results in space forms including roundness for sufficiently pinched initial hypersurfaces evolving by high powers of $F$ in hyperbolic space.)  In this article we extend Li and Lv's work in Euclidean space to the other previously-considered cases of $F$ for corresponding natural conditions on $\Phi$ and we show with sufficient initial curvature pinching evolving hypersurfaces converge in finite time to points that are asymptotically spherical, where the convergence is exponential in the $C^\infty$-topology. 

Let $M_{0}$ be a compact, strictly convex hypersurface of dimension $n\geq 2$, without boundary, smoothly embedded in $\mathbb{R}^{n+1}$ and represented  by some diffeomorphism $X_{0}:\mathbb{S}^{n} \rightarrow 
    X_{0}\left( \mathbb{S}^n \right)= M_{0}\subset \mathbb{R}^{n+1}$.  We 
consider the family of maps $X_{t}=X \left( \cdot, t \right)$ evolving according to 
\begin{equation} \label{E:theflow}
\begin{split}
  \frac{\partial}{\partial t}X \left( x,t \right) 
= - \Phi\left( F\left( \mathcal{W} \left( x,t\right) \right) \right) &\nu
  \left( x,t \right) \mbox{ } x\in \mathbb{S}^n,\mbox{ }  0 < t\leq T\leq \infty \\
X \left( \cdot, 0\right) 
&= X_{0} \mbox{,}
\end{split}
\end{equation}
where $\mathcal{W}\left( x, t\right) $ is the matrix of the Weingarten map of $M_{t}=X_{t}\left( \mathbb{S}^n \right)$ at the point $X_{t}\left( x\right)$ and $\nu\left( x, t\right)$ is the outer unit normal to $M_{t}$ at $X_{t}\left( x\right)$.  Now we describe the properties of the functions $F$ and $\Phi$.

Let $\Gamma_+$ denote the positive cone, $\Gamma_+=\left\{ \kappa = \left( \kappa_1, \ldots, \kappa_n \right) : \kappa_i >0 \mbox{ for all }i=1, \ldots, n \right\}$.  The function $F$ should have the following properties:

\begin{Fconditions} \label{T:Fconds}
\mbox{}
\begin{enumerate}
  \item[\textnormal{a)}] $F\left( \mathcal{W} \right) = f\left( \kappa \left( \mathcal{W} \right) \right)$ where $\kappa\left( \mathcal{W}\right)$ gives the eigenvalues of $\mathcal{W}$ and $f$ is a smooth, symmetric function defined on an open, symmetric cone $\Gamma \subseteq \Gamma_+$. 
  \item[\textnormal{b)}] $f$ is strictly increasing in each argument: $\frac{\partial f}{\partial \kappa_{i}} > 0$ for each $i=1, \ldots, n$ at every point of $\Gamma$.
  \item[\textnormal{c)}] $f$ is homogeneous of degree one: $f\left( k \kappa \right) = k f\left( \kappa\right)$ for any $k>0$.
  \item[\textnormal{d)}] $f$ is strictly positive on $\Gamma$ and $f\left( 1, \ldots, 1\right) = 1$.
  \item[\textnormal{e)}] Either:
    \begin{itemize}
    \item[\textnormal{i)}] $n=2$, or 
      \item[\textnormal{ii)}] $f$ is strictly convex in nonradial directions, or
      \item[\textnormal{iii)}] $f$ is strictly concave in nonradial directions and one of the following hold 
        \begin{enumerate}
            \item[\textnormal{a)}] $f$ approaches zero on the boundary of $\Gamma$,
          \item[\textnormal{b)}] $\sup_{M_{0}} \left( \frac{H}{F} \right) < \liminf_{\kappa \rightarrow \partial \Gamma} \left( \frac{\sum_{i} \kappa_{i}}{f\left( \kappa\right)}\right)$,  
                       \end{enumerate}
         \item[\textnormal{iv)}] $f$ is inverse concave and either
           \begin{itemize}
             \item[\textnormal{a)}] $f_* \rightarrow 0 \mbox{ as } r \rightarrow \partial \Gamma$, or 
            \item[\textnormal{b)}] $\sup_{\substack{\omega \in T_z \mathbb{S}^n \\ \left| \omega\right|=1}} \frac{r\left( \omega, \omega\right) \left( z, 0 \right)}{f_*\left( r\left( z, 0\right) \right)} < \liminf_{r\rightarrow \partial \Gamma} \frac{r_{\max}}{f_*\left( r\right)} \mbox{.}$
            \end{itemize}
          \item[\textnormal{v)}] $f$ satisfies no second derivative condition 
          but either
            \begin{itemize}
              \item[\textnormal{a)}] $M_0$ is axially symmetric, or
              \item[\textnormal{b)}] $M_0$ satisfies a pinching condition of the form
         $$\left| A^0 \right|^2 \leq \sigma H^2 \mbox{,}$$
         where $\sigma$ depends upon $n$ and the second derivative bound on the preserved pinching cone.
         \end{itemize}
 \end{itemize}
\end{enumerate}
\end{Fconditions}

In the above statement we have used $H$ to denote the mean curvature of hypersurface $M_t$; $\left| A^0\right|^2$ is the norm of the trace-free second fundamental form,
$$\left| A^0 \right|^2 = \left| A \right|^2 - \frac{1}{n} H^2 \mbox{,}$$
where $\left| A \right|^2 = \kappa_1^2 + \ldots + \kappa_n^2$ is the norm of the second fundamental form.  We also used $r$ to denote the inverse Weingarten map of $M_t$.  We elaborate on further notation in Section \ref{S:FIE}.

In the case $\Phi\left( F \right) = F$, the speeds $F$ satisfying Conditions \ref{T:Fconds} have all been considered before.  Specifically, condition e) i) in \cite{Asurf}, condition e) ii), e) iii) a) and e) iii) b) in \cite{A94}, e) iv) in \cite{AMZ11}, e) v) a) in \cite{MMW} and condition e) v) b) in \cite{AM12}.  It is not actually necessary to assume $f$ is strictly positive on the positive cone; this follows from property b) above and the Euler identity for homogeneous functions.  As $f$ is degree-one homogeneous, its Hessian is degenerate in the radial direction.\\  

\noindent \textbf{Remark:} We have not included the case of $f$ concave and inverse concave, with no further condition, in part e) above, because the obvious generalisation to the argument in \cite{Apinching} appears to require opposite conditions on $\Phi$ to handle the gradient term and the zero order term in the relevant evolution equation.  Finding an alternative approach to show asymptotic roundness in this case is an open question (\cite{LL20} shows convergence to a point).\\

For our results we will require the function $\Phi: \left[ 0,\infty \right)\rightarrow \mathbb{R}$ to be at least twice differentiable but for convenience and higher regularity of solutions to \eqref{E:theflow} we will assume $\Phi$ is smooth.  For our various results we will require an assortment of the following conditions.  Our requirements do not rule out the possibility that results might be possible with different requirements on $\Phi$ using different approaches.

\begin{Phiconditions} \label{T:Phiconditions}
\begin{itemize}
  \item[\textnormal{a)}] $\Phi\left( 0 \right) = 0$;
  \item[\textnormal{b)}] $\lim_{s\rightarrow \infty} \Phi\left( s\right) = \infty$;
  \item[\textnormal{c)}] $\Phi'\left( s\right) > 0$ for all $s>0$; 
  \item[\textnormal{d)}] Either
  \begin{itemize}
    \item[\textnormal{i)}] $\Phi' s - \Phi \geq 0$ for all $s>0$; or 
    \item[\textnormal{ii)}]  $\Phi's -\Phi \geq \varepsilon \, \Phi' s$ for some $\varepsilon>0$; or 
  \item[\textnormal{iii)}] $\Phi' s - \Phi \leq 0$ for all $s>0$; 
  \end{itemize}
  \item[\textnormal{e)}] $\lim_{s\rightarrow 0} \frac{\Phi'\left( s\right) s^2}{\Phi\left( s\right)} = 0$;
  \item[\textnormal{f)}] $\lim_{s\rightarrow \infty} \frac{\Phi'\left( s\right) s^2}{\Phi\left( s\right)} = \infty$; 
  \item[\textnormal{g)}] $\Phi''\left( s\right) \geq 0$ for all $s>0$, 
  \item[\textnormal{h)}] There exists a constant $c>0$ such that, for all $s$ we have $s \left| \Phi''\left( s\right) \right| \leq c\, \Phi'\left( s\right)$; 
    \item[\textnormal{i)}] $\Phi'' \Phi^2 s + 2 \Phi' \left[ \Phi^2 - \left( \Phi'\right)^2 s^2 \right] > 0$. 
\end{itemize}
\end{Phiconditions}

In particular, Condition c) is needed for parabolicity of \eqref{E:theflow}; this is an essential requirement.   This ensures existence for a short time of a solution to \eqref{E:theflow} (see \cite{B} for example for a suitable short-time existence result).  Condition f) is needed for an upper speed bound via the technique of Chou introduced in \cite{Tso}.  For Condition d), $\Phi' s - \Phi \geq 0$ is used here for the cases $n=2$, for $M_t$ axially symmetric, for $F$ convex and for $F$ concave, while $\Phi' s - \Phi \leq 0$ is used for $F$ inverse concave and either of Conditions \ref{T:Phiconditions} iv) a) or b).  Condition g) above is needed for the case $n=2$ or axially symmetric hypersurfaces; note that Condition g) and a) imply Condition d) i).  Condition h) is needed to absorb $\Phi''$ terms by other good gradient terms in the cases of $F$ strictly concave and strong curvature pinching (Conditions \ref{T:Fconds} e) iii) a) and b) and e) v) b).  Condition i) is used in the case $F$ inverse concave.  It can be checked by calculation that Condition i) is not actually the same as the requirement that $\Phi\circ f$ itself be inverse concave (which we do not require in this article).\\[8pt]

\noindent \textbf{Remarks:}
\begin{enumerate}
  \item Some of the above conditions are similar to those used in \cites{BP, BS, LL20} where constrained curvature flows were considered.  Some of the conditions are new, but arise naturally when we generalise earlier arguments to the case of flow \eqref{E:theflow}.  
    \item The conditions can be slightly weakened if we have other restrictions on our curvature pinching cone of interest; in such cases the conditions for $F$ and for $\Phi$ need only hold on the corresponding domains.
\end{enumerate}

Our main result may now be stated as follows:

\begin{main} \label{T:main}
Let $M_0$ be a smooth, closed, strictly convex, $n$-dimensional hypersurface without boundary, $n\geq 2$, smoothly embedded in $\mathbb{R}^{n+1}$ by $X_0: \mathbb{S}^n \rightarrow \mathbb{R}^{n+1}$.  Let $F$ be a function satisfying Conditions \ref{T:Fconds} and let $\Phi$ be a function satisfying Conditions \ref{T:Phiconditions} a), b), c), e) and f).  Suppose $M_0$ is sufficiently pointwise curvature pinched in the sense that for all $i, j$, 
$$\kappa_i \geq \varepsilon \kappa_j$$
for some $0< \varepsilon < 1$ depending on $F$ and on $\Phi$.   Then, if (at least) one of the following hold:
\begin{itemize}
  \item In the cases of Conditions \ref{T:Fconds} e) i) and e) v) a), if $\Phi$ also satisfies Conditions \ref{T:Phiconditions} g);
  \item In the cases of Conditions \ref{T:Fconds} e) ii) and e) iii) a) and b), if $\Phi$ also satisfies Conditions \ref{T:Phiconditions} d) i) and h);
\item In the case of Conditions \ref{T:Fconds} e) iv) a) and b), if $\Phi$ also satisfies Conditions \ref{T:Phiconditions} d) iii) and Conditions \ref{T:Phiconditions} i);
\item In the case of Conditions \ref{T:Fconds} e) v) b), if $\Phi$ also satisfies Conditions \ref{T:Phiconditions} d) ii) and h);
\end{itemize}
 there exists a family of smooth, strictly convex hypersurfaces $\left\{ M_t = X_t\left( \mathbb{S}^n \right) \right\}_{0\leq t < T < \infty}$, unique up to reparametrisation, satisfying \eqref{E:theflow}, with initial condition $X\left( x, 0 \right) = X_0\left( x\right)$.  The image hypersurfaces $M_t$ contract to a point $p\in \mathbb{R}^{n+1}$ at time $T$.  Under appropriate rescaling the image hypersurfaces smoothly converge uniformly and exponentially to spheres as the rescaled time parameter $\tau \rightarrow \infty$. 
\end{main}

\noindent \textbf{Remark:} In the case of Conditions \ref{T:Fconds} e) i) above, the curvature pinching requirement may be stated more explicitly as
 $$\frac{\kappa_{\max}}{\kappa_{\min}} \leq 1 + \frac{2\Phi'\left( f\left( \kappa\right) \right)}{\Phi''\left( f\left( \kappa\right) \right)\, f\left( \kappa\right)} $$
while in the case of Conditions \ref{T:Fconds} e) v) a) the requirement is
$$\frac{\kappa_{\mbox{axial}}}{\kappa_{\mbox{rotational}} }\leq 1 + \frac{2\Phi'\left( f\left( \kappa\right) \right)}{\Phi''\left( f\left( \kappa\right) \right)\, f\left( \kappa\right)} \mbox{,}$$
where $\kappa_{\mbox{axial}}$ and $\kappa_{\mbox{rotational}}$ denote the axial and rotational curvatures respectively of the axially-symmetric hypersurface $M_0$.\\[8pt]

We will use similar notation as in \cites{Apinching, M05, AMZ11} and elsewhere.  The structure of this article is as follows.  In Section \ref{S:examplespeeds} we give some example flow speeds satisfying our conditions.  In Section \ref{S:FIE} we detail some flow independent geometric results that will be needed in the analysis.  In Section \ref{S:evlneqns} we list the evolution equations we will need in the subsequent analysis.  These are all derived in a straightforward manner as in previous work.  In Section \ref{S:pinching} we establish the crucial curvature pinching estimate for each of the cases of $F$ from Conditions \ref{T:Fconds} part e) in turn.  In Section \ref{S:T} we discuss regularity and complete the proof of the main theorem.

In a forthcoming article \cite{M20} we will consider the case of hypersurfaces that contract self-similarly under general nonhomogeneous curvature functions.  By our approach, to characterise such contracting hypersurfaces as spheres still requires curvature pinching but generally requires fewer structure conditions on $\Phi$ than Conditions \ref{T:Phiconditions}.  For particular $F$ and $\Phi$, generally involving the Gauss curvature, improved results that do not require curvature pinching are available.

\section{Example speeds} \label{S:examplespeeds}
\newtheorem{caic}{Lemma}[section]

The class of $F$ considered in this article is by now quite well known.  Many such $F$ are given for example in \cites{A94, Apinching, AMZ11}.  The functions $\Phi$ are less well-known, but some have been considered before in \cites{BS, BP, LL20}.  Some suitable functions $\Phi$ are:

\begin{enumerate}
  \item $\Phi\left( s\right) = \sum_{i=1}^{\ell} c_i s^{k_i}$, for some $c_i>0$, $k_i>0$.  We assume $\ell \in \mathbb{N}\backslash \left\{ 1\right\}$ as the case $\ell=1$ has been considered elsewhere.  Clearly these $\Phi$ satisfy Conditions \ref{T:Phiconditions} a), b) and c).  We have 
  $$\Phi'\left( s\right) = \sum_{i=1}^{\ell}  c_i k_i\, s^{k_i-1}$$
  so	
  $$\Phi'\left( s\right) s - \Phi\left( s\right) = \sum_{i=1}^{\ell} c_i \left( k_i-1 \right) s^{k_i} \mbox{.}$$
  Thus for Condition \ref{T:Phiconditions} d) we will have 
  \begin{itemize}
    \item $\Phi'\, s - \Phi \geq 0$ for all $s>0$ if all $k_i \geq 1$;
    \item  $\Phi'\, s - \Phi \geq \varepsilon \, \Phi'\, s$ for all $s>0$ if all $k_i \geq \frac{1}{1-\varepsilon}$;
    \item $\Phi'\, s - \Phi \leq 0$ for all $s>0$ if all $k_i$ satisfy $0<k_i \leq 1$.  
  \end{itemize}
  
  Conditions e) and f) are clearly satisfied by these $\Phi$.  For Condition \ref{T:Phiconditions} g), we will have $\Phi''\left( s\right) \geq 0$ for all $s>0$ if $k_i \geq 1$ for all $i$.  Condition h) is satisfied by these $\Phi$ for any $c\geq \max_i c_i \left| k_i -1\right|$.  Condition i) is satisfied if $k_i>1$ for all $i$.\\
   
  \item $\Phi\left( s\right) = \ln\left( 1 + s\right)$.  Clearly $\Phi$ satisfies Conditions \ref{T:Phiconditions} a), b) and c).  For d) we have $\Phi'\, s - \Phi <0$.  Conditions e) and f) are again clearly satisfied.  Condition h) is also satisfied for any $c\geq 1$ and Condition i) is also satisfied.\\
  
   \item $\Phi\left( s\right) = e^s -1$.  Clearly $\Phi$ satisfies Conditions \ref{T:Phiconditions} a), b) and c).  For d) we have $\Phi'\, s - \Phi \geq 0$.  Conditions e) and f) are again clearly satisfied while for g) we have $\Phi''\left( s\right) >0$ for all $s$.  Neither h) nor i) are satisfied for this $\Phi$.\\
   
   \item More $\Phi$ can be built from linear combinations of the above.  For example, $\Phi\left( s\right)= \ln\left( 1 + s\right) + s^2$ clearly satisfies Conditions \ref{T:Phiconditions} a), b), c), d) i), e), f), g) and h).\\
  
  \item In applications, other $\Phi$ might be constructed piecewise; of course the regularity of solutions to \eqref{E:theflow} will be affected (for non-smooth degree-one homogeneous speeds we refer the interested reader to \cites{AHMWWW14, AM16} for the cases of convex speeds and $n=2$ respectively).\\
\end{enumerate}

\section{Notation, time independent facts and geometric estimates} \label{S:FIE}

\newtheorem{AndW}{Theorem}[section]
\newtheorem{Consign}[AndW]{Corollary}
\newtheorem{Strictcon}[AndW]{Corollary}
\newtheorem{Fconvprops}[AndW]{Corollary}
\newtheorem{Fprops}[AndW]{Lemma}
\newtheorem{cop}[AndW]{Theorem}
\newtheorem{rin}[AndW]{Corollary}
\newtheorem{HCinequal}[AndW]{Lemma}
\newtheorem{gradA}[AndW]{Lemma}
\newtheorem{tight}[AndW]{Lemma}
\newtheorem{trace}[AndW]{Lemma}
\newtheorem{extreme}[AndW]{Lemma}
\newtheorem{AMgeom}[AndW]{Lemma}
\newtheorem{AMHa}[AndW]{Lemma}
\newtheorem{conseq}[AndW]{Corollary}
\newtheorem{Tso}[AndW]{Corollary}
\newtheorem{hbound}[AndW]{Corollary}
\newtheorem{ratio}[AndW]{Lemma}
\newtheorem{radii}[AndW]{Corollary}

Our notation is the same as that used elsewhere (egs \cites{Asurf, AMZ11, Hu}).  In particular we write $g = \left\{ g_{ij}\right\}$, $A = \left\{ h_{ij}\right\}$ and $\mathcal{W} = \left\{ h^{i}_{\> j} \right\}$ to denote respectively the metric, second fundamental form and Weingarten map of the evolving $M_t$.  The mean curvature of $M_t$ is
$$H= g^{ij}h_{ij} = h^{i}_{\> i}$$
and the norm of the second fundamental form is
$$\left| A \right|^{2} = g^{ij}g^{lm}h_{il}h_{jm} = h^{j}_{\> l}h_{j}^{\;\, l}$$
where $g^{ij}$ is the $\left( i, j\right)$-entry of the inverse of the matrix $\left( g_{ij}\right)$.  As mentioned earlier, the norm of the trace-free component of the second fundamental form,
$$\left| A^{0}\right|^{2} = \left| A \right|^{2} - \frac{1}{n} H^{2} = \frac{1}{n} \sum_{i<j} \left( \kappa_{i} - \kappa_{j} \right)^{2} \mbox{,}$$
is identically equal to zero when $M_t$ is a sphere.  As in earlier work we will also set 
$$C = \kappa_{1}^{3} + \ldots + \kappa_{n}^{3} \mbox{,}$$
where the letter $C$ is chosen simply by convention.  Throughout this paper we sum over repeated indices from $1$ to $n$ unless otherwise indicated.  Raised indices indicate contraction with the metric.

We denote by $\left( \dot{F}^{kl} \right)$ the matrix of first partial derivatives of $F$ with respect to the components of its argument:
$$\left. \frac{\partial}{\partial s} F\left( A+sB\right) \right|_{s=0} = \dot F^{kl} \left( A\right) B_{kl} \mbox{.}$$
Similarly for the second partial derivatives of $F$ we write
$$\left. \frac{\partial^{2}}{\partial s^{2}} F\left( A+sB\right) \right|_{s=0} = \ddot F^{kl, rs} \left( A\right) B_{kl} B_{rs} \mbox{.}$$
Throughout the article unless the argument is explicitly indicated we will always evaluate partial derivatives of $F$ at $\mathcal{W}$ and partial derivatives of $f$ at $\kappa\left( \mathcal{W}\right)$.  We will further use the shortened notation $\dot f^{i} = \frac{\partial f}{\partial \kappa_{i}}$ and $\ddot f^{ij} = \frac{\partial^{2} f}{\partial \kappa_{i} \partial \kappa_{j}}$ where appropriate.

Let us next review some results for symmetric functions of the eigenvalues of matrices.  The first theorem, proofs of which appear in \cite{AAniso} and in \cite{Ge2}, allows us to swap between derivatives of $F$ and derivatives of $f$ at a particular point (see also \cite{A94} equation (2.23)).  We denote by $\mbox{Sym}\left( n\right)$ the set of all $n\times n$ real symmetric matrices.

\begin{AndW} \label{T:AndW}
Let $f$ be a $C^{2}$ symmetric function defined on a symmetric region $\Omega$ in $\mathbb{R}^{n}$.  Let $\tilde \Omega = \left\{ A\in \mbox{\textnormal{Sym}}\left(n \right): \kappa \left( A\right) \in \Omega \right\}$ and define $F: \tilde \Omega \rightarrow \mathbb{R}$ by $F\left( A\right) = f\left( \kappa\left( A\right)\right)$.  Then at any diagonal $A\in \tilde \Omega$ with distinct eigenvalues, the second derivative of $F$ in direction $B\in \mbox{\textnormal{Sym}}\left( n\right)$ is given by
$$\ddot F^{kl, rs} B_{kl}B_{rs} = \sum_{k,l} \ddot f^{kl} B_{kk}B_{ll} + 2\sum_{k<l} \frac{\dot f_{k} - \dot f_{l}}{\kappa_{k} - \kappa_{l}} B_{kl}^{2} \mbox{.}$$
\end{AndW}

From the above result we have immediately:

\begin{Consign} 
  $F$ is convex (concave) at $\mathcal{W}$ if and only if $f$ is convex (concave) at $\kappa\left( \mathcal{W}\right)$ and
$$\frac{\dot f_{i} - \dot f_{j}}{\kappa_{i} - \kappa_{j}} \geq \left( \leq\right) 0 \mbox{ for all } i \ne j \mbox{.}$$
\end{Consign}

We also have the following result which is useful when the convexity or concavity of $F$ is strict in nonradial directions.

\begin{Consign}[Lemma 7.12, \cite{A94}] \label{T:consign}
If $F$ satisfies Conditions \ref{T:Fconds} and is strictly convex in nonradial directions, then there exists $\underline{C}>0$ such that
$$\ddot F^{kl,rs} B_{kl} B_{rs} \geq \underline{C} \frac{\left| \nabla A\right|^2}{\left| A \right|} \mbox{.}$$
On the other hand, if $F$ satisfies Conditions \ref{T:Fconds} and is strictly concave in nonradial directions, then there exists $\overline{C} > 0$ such that
$$\ddot F^{kl,rs} B_{kl} B_{rs} \leq -\underline{C} \frac{\left| \nabla A\right|^2}{\left| A \right|} \mbox{.}$$
\end{Consign}

The inequalities of the following Lemma are derived in \cite{U} for concave $F$.  The signs are opposite for convex $F$.

\begin{Fprops} \label{T:Fprops}
For any convex (concave) $F$ satisfying Conditions \ref{T:Fconds}, for all $\kappa \in \Gamma$,
\begin{itemize}
  \item[\textnormal{i)}] $f\left( \kappa \right) \geq \left( \leq \right) \frac{1}{n} H$,
  \item[\textnormal{ii)}] $\sum_{k} \dot{f}_{k} = \mbox{\textnormal{trace}}\left( \dot{F}^{kl} \right) \leq \left( \geq \right) 1$. 
\end{itemize}
\end{Fprops}

The next fact follows from the structure conditions on $F$.  It is shown for example in \cite{U}.

  \begin{trace} \label{T:trace}
  For any concave function $F$ satisfying Conditions \ref{T:Fconds} a) to d),
  $$\sum_{k=1}^n \dot f_k = \mbox{\textnormal{trace}} \left( \dot F^{kl} \right) \geq 1 \mbox{.}$$
  \end{trace}

A proof of the next result appears in \cite{M05}.

\begin{Fconvprops} \label{T:Fconvprops}
\mbox{}
\begin{itemize}
  \item[\textnormal{i)}] If $F$ satisfies Conditions \ref{T:Fconds} and $F$ is convex (concave) at $\mathcal{W}$, then at this $\mathcal{W}$,
  $$F \left| A\right|^{2} - H \dot{F}^{kl} h_{km} h^{m}_{\ \> l} \leq \left( \geq \right) 0 \mbox{.}$$
  \item[\textnormal{ii)}] If $F$ satisfies Conditions \ref{T:Fconds} and $F$ is convex (concave) at $\mathcal{W}$, then at this $\mathcal{W}$,
  $$FC - \left| A\right|^{2} \dot{F}^{kl} h_{km} h^{m}_{\ \> l}  \leq \left( \geq \right) 0 \mbox{.}$$
\end{itemize}
\end{Fconvprops}

 The next inequalities appear as \cite[Lemma 2.3 (i)]{Hu}, \cite[Lemma 1.4 (iii)]{HV} 
and Lemma \cite[Lemma 4.1]{CRS} respectively.  The third is attributed to Huisken and is a stronger version of an inequality in \cite{Hu}.  It is important here because it contains the full $\nabla A$ on the right hand side rather than $\left| \nabla H\right|$.\\

\begin{HCinequal} \label{T:HCinequal} 
  If $H>0$ and $h_{ij} \geq \varepsilon H g_{ij}$ is valid for some $\varepsilon >0$ then, using coordinates at any particular point that diagonalise the Weingarten map, we have
\begin{itemize}
  \item[\textnormal{(i)}] $HC - \left( \left| A\right|^2 \right)^2 = \sum_{i<j}^n \kappa_i \kappa_j \left( \kappa_i - \kappa_j \right)^2 \geq n\, \varepsilon^2 H^2 \left| A^0\right|^2  \mbox{;}$
  \item[\textnormal{(ii)}] $n\, C - H \left| A\right|^2 = \frac{1}{2} \sum_{i\ne j} \left( \kappa_i + \kappa_j \right) \left( \kappa_i - \kappa_j \right)^2 \geq 2\,n\,\varepsilon \, H \left| A^0 \right|^2 \mbox{;}$
  \item[\textnormal{(iii)}] $\left| H \nabla_i h_{jk} - h_{jk} \nabla_i H \right|^2 \geq \frac{n-1}{2} \varepsilon^2 H^2 \left| \nabla A\right|^2 \mbox{,}$
\end{itemize}
  where $\left| A \right|^2 = \kappa_1^2 + \ldots + \kappa_n^2$ is the squared norm of the second fundamental form, $\left| A^0\right|= \left| A\right|^2 - \frac{1}{n} H^2$ is the squared norm of the trace-free second fundamental form and $C= \kappa_1^3 + \ldots + \kappa_n^3$.
  \end{HCinequal}
  
We also have from \cite{Hu},
\begin{gradA} \label{T:gradA}
$$ \left| \nabla A \right|^2 \geq \frac{3}{n+2} \left| \nabla H \right|^2$$
and equivalently
$$\left| \nabla A^0 \right|^2 \geq \frac{2\left( n-1\right)}{3n} \left| \nabla A\right|^2 \mbox{.}$$
\end{gradA}

The following results, that holds for tight curvature pinching, are \cite[Lemma 2.2, Lemma 2.3]{AM12}.  They will be used in the case of Conditions \ref{T:Fconds} e) v) b).  Since $\Gamma$ is open, there exists $\sigma \in \left( 0, \frac{1}{n\left( n-1\right)} \right)$ such that
$$\Gamma_0 = \left\{ \kappa\left( \mathcal{W}\right) : \left| A^0 \right|^2 \leq \sigma_0 \, H^2 \right\} \subset \Gamma \mbox{.}$$

\begin{tight} \label{T:tight}
If $H>0$ and $\left| A^0\right|^2 = \varepsilon H^2$ at $p$ with $\varepsilon < \frac{1}{n\left( n-1\right)}$, then the second fundamental form is positive definite at $p$ and we have
\begin{itemize}
  \item[\textnormal{(i)}] $\left( 1 - \sqrt{n\left( n-1\right) \varepsilon} \right) \frac{H}{n} \leq \kappa_i \leq \left( 1 + \sqrt{n\left( n-1\right) \varepsilon} \right) \frac{H}{n}$;
  \item[\textnormal{(ii)}] $n\, C - \left( 1 + n\, \varepsilon \right) H \left| A\right|^2 \geq \varepsilon \left( 1 + n\, \varepsilon \right) \left( 1 - \sqrt{n\left( n-1\right) \varepsilon} \right) H^3$.\\[8pt]
  \end{itemize}
\end{tight}

\noindent \textbf{Remarks:} 
\begin{enumerate}
  \item As a consequence of Lemma \ref{T:tight}, (i), the curvatures of $M$ are pinched at $p$ in the sense that
$$\kappa_{\min} \geq \frac{\left( 1 - \sqrt{n\left( n-1\right) \sigma} \right)}{\left( 1 + \sqrt{n\left( n-1\right) \sigma} \right)} \kappa_{\max} \mbox{.}$$
  \item  For the inequality $h_{ij} \geq \varepsilon H g_{ij}$ to hold on a convex hypersurface, it follows by taking the trace that $\varepsilon \leq \frac{1}{n}$.  Inequality (i) of Lemma \ref{T:tight} says that under the same assumptions, we have 
\begin{equation} \label{E:pinchequiv}
  h_{ij} \geq \varepsilon H g_{ij} \mbox{ with } \varepsilon = \frac{1 - \sqrt{n\left( n-1\right) \sigma}}{n} \mbox{.}
\end{equation} 
\end{enumerate}

For the case of Conditions \ref{T:Fconds} e) v) b) we will also need the following estimates on $F$ that follow via integration as in \cite{AM12}, Section 4.  

\begin{AMHa} \label{T:AMHa}
For arbitrary $\kappa\left( \mathcal{W} \right) \in \Gamma_0$ and arbitrary symmetric $2$-tensors $P$ and $Q$, there exists a $\mu\geq 0$ such that
$$\left| \ddot F^{kl,rs} \left( \mathcal{W} \right) P_{kl} Q_{rs} \right| \leq \frac{\mu}{H} \left| P \right| \left| Q\right| \mbox{.}$$
Moreover, we have the matrix inequalities
$$\left( 1 - \mu \frac{\left| A^0\right|}{H} \right) I \leq \dot F \left( \mathcal{W}\right) \leq \left( 1 + \mu \frac{\left| A^0\right|}{H} \right) I$$
and the estimate
$$\left| F\left( \mathcal{W}\right) - H \right| \leq \frac{\mu}{2} \frac{\left| A^0\right|^2}{H} \mbox{.}$$
\end{AMHa}

\noindent \textbf{Remark:} The constant $\mu\geq 0$ is determined by the form of $F$.  Given $\mu$, with a pinching requirement of the form $\left| A^0 \right|^2 \leq \sigma H^2$, if we take $\sqrt{\sigma} < \frac{1}{\mu}$ we can be sure that $\dot F$ is positive definite.\\

The last result we state in this section will be needed for the case of Conditions \ref{T:Fconds} e) v) a).  It is proved as \cite[Lemma 5.1]{MMW}.

\begin{extreme} \label{T:extreme}
   Let $G\left( \mathcal{W} \right) = g\left( \kappa\left( \mathcal{W} \right) \right)$ be a smooth, symmetric, homogeneous of degree zero function in the principal curvatures of an axially symmetric hypersurface, where the coordinates are chosen such that $x_1$ is the axial direction.  At any stationary point of $G$ for which $\dot G$ is nondegenerate,
  $$\left( \dot G^{ij} \ddot F^{kl,rs} - \dot F^{ij} \ddot G^{kl,rs} \right) \nabla_i h_{kl} \nabla_{j} h_{rs} = \frac{2\, f\, \dot g^1}{\kappa_2 \left( \kappa_2 - \kappa_1 \right)} \left( \nabla_1 h_{22} \right)^2 \mbox{.}$$
\end{extreme}

\noindent \emph{Note:} By nondegenerate we mean that in coordinates that diagonalise the Weingarten map, the corresponding diagonal matrix of $\dot G$ has no zero entries on the diagonal.

\section{Evolution equations} \label{S:evlneqns}
\newtheorem{evlneqns}{Lemma}[section]
\newtheorem{evlneqns2}[evlneqns]{Lemma}

The following evolution equations may be derived similarly as in previous work (see, for example \cites{Hu, A94, LL20}).  We write $u\left( x, t\right) = \left< X\left( x, t\right), \nu\left( x, t\right) \right>$ for the support function of the evolving hypersurface $M_t$.  We also set $Z= \frac{\Phi}{u- \delta}$ for constant $\delta$ to be chosen, $Z_\sigma = \left| A^0\right|^2 - \sigma H^2$ and denote by $G$ a general degree-$\alpha$ homogeneous symmetric function of the Weingarten map. 

\begin{evlneqns} \label{T:evlneqns}
Under the flow \eqref{E:theflow} we have the following evolution equations for geometric quantities associated with the evolving hypersurface $M_t$:
\begin{itemize}
  \item[\textnormal{a)}] $\frac{\partial}{\partial t} g_{ij} = -2 \Phi  \, h_{ij}$; $\frac{\partial}{\partial t} g^{ij} = 2 \Phi\, h^{ij}$;
  \item[\textnormal{b)}] $\frac{\partial}{\partial t} h_{ij} = \Delta_{\dot \Phi} h_{ij} + \ddot \Phi^{kl, rs} \nabla_i h_{kl} \nabla_j h_{rs} + \dot \Phi^{kl} k_{k}^{\ p} h_{pl} h_{ij} -\left( \Phi' F + \Phi\right) h_{i}^{\ k}h_{kj}$;
  \item[\textnormal{c)}] $\frac{\partial}{\partial t} h^i_{\ j} = \Delta_{\dot \Phi} h^i_{\ j} + \ddot \Phi^{kl, rs} \nabla^i h_{kl} \nabla_j h_{rs} + \dot \Phi^{kl} k_{k}^{\ p} h_{pl} h^{i}_{\ j} -\left( \Phi' F - \Phi\right) h^{ik}h_{kj}$;
  \item[\textnormal{d)}] $\frac{\partial}{\partial t} G = \Delta_{\dot \Phi} G + \left( \dot G^{ij} \ddot \Phi^{kl, rs} - \dot \Phi^{ij} \ddot G^{kl, rs} \right) \nabla_i h_{kl} \nabla_{j} h_{rs} + \alpha \, \dot \Phi^{kl} h_k^{\ p}h_{pl} G$
  
   \hspace*{0.8cm}$- \left( \Phi' F - \Phi \right) \dot G^{ij} h_i^{\ k} h_{kj}$;
  \item[\textnormal{e)}] $\frac{\partial}{\partial t} F = \Delta_{\dot \Phi} F + \Phi'' \left| \nabla F \right|_{\dot F}^2 + \dot F^{kl} h_k^{\ p}h_{pl} \Phi$;
  \item[\textnormal{f)}] $\frac{\partial}{\partial t} H = \Delta_{\dot \Phi} H + \ddot \Phi^{kl, rs} \nabla^i h_{kl} \nabla_i h_{rs} + \dot \Phi^{kl} h_k^{\ p}h_{pl} H  -\left( \Phi' F - \Phi\right) \left| A\right|^2$;
  \item[\textnormal{g)}] $\frac{\partial}{\partial t} \Phi = \Delta_{\dot \Phi} \Phi + \Phi' \dot F^{kl} h_k^{\ p}h_{pl} \Phi$;
  \item[\textnormal{h)}] $\frac{\partial}{\partial t} \left( \frac{H}{F} \right) = \Delta_{\dot \Phi} \left( \frac{H}{F} \right) + \frac{2}{F} \dot \Phi^{kl} \nabla_k F \nabla_l \left( \frac{H}{F} \right) + \frac{\Phi'}{F} \ddot F^{kl, rs} \nabla^i h_{kl} \nabla_i h_{rs}$
  
  \hspace*{0.8cm} $+ \frac{\Phi''}{F^2} \left( F g^{ij} - H \dot F^{ij} \right) \nabla_i F \nabla_j F + \frac{1}{F^2} \left( \Phi' F - \Phi \right) \left( H \dot F^{kl} h_k^{\ p} h_{pl} - F\left| A\right|^2 \right)$;
  \item[\textnormal{i)}] $\frac{\partial}{\partial t} \left( \frac{K}{F^n} \right) = \Delta_{\dot \Phi} \left( \frac{K}{F^n} \right) + \Phi' \dot F^{kl} \left( \frac{n+2}{F}\nabla_k F - \frac{1}{K} \nabla_k K \right)  \nabla_l  \left( \frac{K}{F^n} \right)$
  
  \hspace*{0.8cm}    $+ \Phi' \frac{K}{F^n} \left( h^{-1}\right)^{ij} \ddot F^{kl, rs} \nabla_i h_{kl} \nabla_j h_{rs} + \Phi'' \frac{K}{F^{n+1}}\left[ F \left( h^{-1}\right)^{ij} - n  \dot F^{ij} \right] \nabla_i F \nabla_j F$
  
  \hspace*{0.8cm} $+ \frac{K}{F^{n+2}} \dot \Phi^{kl} \left( h^{-1}\right)_{pq}  \left( h^{-1}\right)_{rs} \left( \Phi \nabla_k h^{pr} - h^{pr} \nabla_k \Phi \right) \left( \Phi \nabla_l h^{qs} - h^{qs} \nabla_l \Phi \right)$
  
  \hspace*{0.8cm} $+ \left( \Phi' F - \Phi \right)\frac{K}{F^{n+1}}  \left( n \dot F^{kl} h_k^{\ p} h_{pl} - H\, F \right)$;
  \item[\textnormal{j)}] $\frac{\partial}{\partial t}  u = \Delta_{\dot \Phi} u + \Dot \Phi^{kl} h_k^{\ p} h_{pl} u - \Phi F' - \Phi$;
  \item[\textnormal{k)}] $\frac{\partial}{\partial t}  Z = \Delta_{\dot \Phi} Z + \frac{2}{u-\delta} \dot \Phi^{kl} \nabla_k \left( u-\delta\right) \nabla_l Z - \frac{\delta}{u-\delta} \Phi' \dot F^{kl} h_k^{\ p} h_{pl} Z + \left( \Phi' F + \Phi \right) \frac{Z}{u-\delta}$;
  \item[\textnormal{l)}] $\frac{\partial}{\partial t}  Z_\sigma = \Delta_{\dot \Phi} Z_\sigma + 2\left[ h^{ij} - \left( \frac{1}{n} + \sigma \right) H g^{ij} \right] \ddot \Phi^{kl, rs} \nabla_i h_{kl} \nabla_j h_{rs}$
  
  \hspace*{0.8cm} $-2\dot \Phi^{ij} \left( \nabla_i h_{kl} \nabla_j h^{kl} - \frac{1}{n} \nabla_i H \nabla_j H \right)+ 2\sigma \dot \Phi^{ij} \nabla_i H \nabla_j H + 2  \dot \Phi^{kl} h_k^{\ p} h_{pl} Z_\sigma$
  
  \hspace*{0.8cm} $+ \frac{2}{n} \left( \Phi - \Phi' F\right) \left[ n \, C - \left( 1 + n\sigma\right) H \left| A\right|^2 \right]$.
  \end{itemize}
  \end{evlneqns}

Above we use the notation as in \cite{LL20}:
$$\dot \Phi^{kl} = \Phi' \dot F^{kl}$$
$$\ddot \Phi^{kl, rs} = \Phi' \ddot F^{kl, rs} + \Phi'' \dot F^{kl} \dot F^{rs}$$
$$\left| T \right|_{\dot F} = \dot F^{kl} T_k T_l$$
$$\Delta_{\dot \Phi} = \dot \Phi^{kl} \nabla_k \nabla_l$$
where $'$ denotes derivative of a function of one variable with respect to its argument.\\

In the case of Conditions \ref{T:Fconds} e) iv), as in earlier work with this condition on $F$, we instead parametrise the evolving hypersurface using its support function, incorporating a tangential diffeomorphism into \eqref{E:theflow} such that the parametrisation is preserved.  Details of this procedure may be found in \cites{AHarnack, AGauss, AMZ11, U} for examples.  Then $M_t$ has support function $s: \mathbb{S}^n \times \left[ 0, T\right) \rightarrow \mathbb{R}^{n+1}$ given by
\begin{equation} \label{E:s}
  s\left( x, t\right) = \left< X\left( x, t\right), x\right> \mbox{,}
\end{equation}
where $x$ is the outer unit normal to $M_t$ at $X\left( x, t\right)$ for all $t\in \left[ 0, T\right)$.  The matrix of the inverse Weingarten map $\mathcal{W}^{-1}$ has entries given by
$$r_{ij} = \overline{\nabla}_i  \overline{\nabla}_j s + \overline{g}_{ij} s \mbox{,}$$
where $\overline{\nabla}$ denotes the covariant derivative on $\mathbb{S}^n$.  The eigenvalues of $\mathcal{W}^{-1}$ are the principal radii of curvature of $M_t$ that we will denote by $r_1, \ldots r_n$.

With the appropriate tangential term added to \eqref{E:theflow} to ensure \eqref{E:s} is preserved, we have the following evolution equations:

\begin{evlneqns2} \label{T:evlneqns2}
\begin{itemize}
  \item[\textnormal{a)}] $\frac{\partial}{\partial t} s= - \Phi \left( F_*^{-1}\left( r_{ij}\right) \right)$;
  \item[\textnormal{b)}] $\frac{\partial}{\partial t} r_{ij} = \Phi' F_*^{-2} \dot F_*^{pq} \overline{\nabla}_p \overline{\nabla}_q r_{ij} - \ddot \Phi^{pq, rs} \overline{\nabla}_i r_{pq} \overline{\nabla}_j r_{rs} - \Phi' F_*^{-2} \mbox{tr } \dot F_* \cdot r_{ij} + \left( \Phi' F_*^{-1} - \Phi \right) \overline{g}_{ij}$;
  \item[\textnormal{c)}] $\frac{\partial}{\partial t} \Phi = \Phi' F_*^{-2} \dot F^{pq} \overline{\nabla}_p \overline{\nabla}_q \Phi + \Phi' F_*^{-2} \mbox{\textnormal{tr} }\dot F_* \, \Phi$;
    \item[\textnormal{d)}] $\frac{\partial}{\partial t} \Phi^{-1} = \Phi' F_*^{-2} \dot F^{pq} \overline{\nabla}_p \overline{\nabla}_q \Phi^{-1} - \frac{2\Phi'}{\Phi^3} F_*^{-2} \dot F_*^{pq} \overline{\nabla}_p \Phi \overline{\nabla}_q \Phi - \frac{\Phi'}{\Phi} F_*^{-2} \mbox{\textnormal{tr} }\dot F_*$.
    \end{itemize}
    \end{evlneqns2}

\section{Pinching estimates} \label{S:pinching}

Curvature pinching of the type $\kappa_i \geq \varepsilon \kappa_j$ may be established in each of the cases of Conditions \ref{T:Fconds} e) using suitable adjustments of previous work.  Natural requirements on $\Phi$ soon become clear in each case.

Consequences of curvature pinching include bounds on degree-zero homogeneous functions of the principle curvatures: in particular, curvature pinching implies existence of constants $0< \underline{C} \leq \overline{C}$ such that $\dot F$ is uniformly elliptic:
$$0< \underline{C} I \leq \dot F \leq \overline{C} I \mbox{.}$$
Also, as in \cite{A94}, curvature pinching implies the ratio of the outer radius to the inradius of $M_t$ is bounded above.

\subsection{Case of Conditions \ref{T:Fconds} e) i) and e) v) a)}
 Both these cases can be handled using the degree-zero homogeneous function
 $$G\left( \mathcal{W}\right) = \frac{n \left| A^0 \right|^2}{H^2}\mbox{.}$$
 For the axially symmetric case, let us denote by $\kappa_1$ the curvature in the axial direction, and $\kappa_2$ the rotational curvature.  In either case we may write 
 $$g\left( \kappa\right) = \frac{\left( n-1\right) \left( \kappa_1 - \kappa_2\right)^2}{\left[ \kappa_1 + \left( n-1\right) \kappa_2 \right]^2} \mbox{.}$$
 It is important to keep in mind that $g$ is symmetric in $\left( \kappa_1, \kappa_2\right)$ in the $n=2$ case, but not in the general axially symmetric case.  We have
 $$\dot g^1 = \frac{2n\left( n-1\right) \kappa_2 \left( \kappa_1 - \kappa_2\right)}{H^3} \mbox{ and } \dot g^2 = \frac{2n\, \kappa_1\left( \kappa_2 - \kappa_1\right)}{H^3} \mbox{.}$$
 
We wish to show that under the flow \eqref{E:theflow}, the maximum of $G$ does not increase.  Since $G$ is degree-zero homogeneous, $\Phi' F - \Phi \geq 0$ and 
$$\dot G^{ij} h_i^{\ k} h_{kj} = \dot g^1 \kappa_1^2 + \left( n-1\right) \dot g^2 \kappa_2^2 = \frac{2n\left( n-1\right) \kappa_1 \kappa_2}{H^3} \left( \kappa_1 -\kappa_2\right)^2 \geq 0 \mbox{,}$$
 the zero order terms in the evolution equation for $G$ are nonpositive.  It remains to check the gradient term.  Since $F$ is degree-one homogeneous and $G$ is degree-zero homogeneous, using the condition $\nabla G = 0$ at a maximum we find using cancellation as in \cite{Asurf} that in the present setting  
\begin{align*}
& \left( \dot G^{ij} \ddot \Phi^{kl, rs} - \dot \Phi^{ij} \ddot G^{kl, rs} \right) \nabla_i h_{kl} \nabla_j h_{rs}\\
  &= \frac{F \dot g^1}{\kappa_1 \kappa_2^2 \left( \kappa_2 - \kappa_1\right)}\left\{ \left[ \Phi'' F \kappa_1\left( \kappa_2 - \kappa_1 \right) + 2 \Phi' \kappa_1 \kappa_2 \right] \left( \nabla_1 h_{22} \right)^2 \right. \\
  & \qquad \left. + \left[ - \frac{\Phi'' F \kappa_2}{\left( n-1\right)} \left( \kappa_2 - \kappa_1\right) + 2 \Phi' \kappa_1 \kappa_2 \right]  \left( \nabla_2 h_{11} \right)^2 \right\}\\
  &=\frac{-2n\left( n-1\right)F}{\kappa_1 \kappa_2 H^3}\left\{ \kappa_1 \left[ \Phi'' F \left( \kappa_2 - \kappa_1 \right) + 2 \Phi'  \kappa_2 \right] \left( \nabla_1 h_{22} \right)^2 \right. \\
  & \qquad \left. + \kappa_2 \left[ - \frac{\Phi'' F }{\left( n-1\right)} \left( \kappa_2 - \kappa_1\right) + 2 \Phi' \kappa_1 \right]  \left( \nabla_2 h_{11} \right)^2 \right\} \mbox{.}
  \end{align*}

 In view of symmetry, in the case $n=2$ we may assume $\kappa_2 \geq \kappa_1$ and for the above gradient term to be negative we require
 $$-2\Phi' \kappa_1 + \Phi'' F \left( \kappa_2 - \kappa_1 \right) \leq 0 \mbox{ and } 2 \Phi' \kappa_2 + \Phi'' F \left( \kappa_2 - \kappa_1 \right) \geq 0 \mbox{.}$$
If it happens that $\Phi''=0$ at the maximum of $G$, then these conditions are obviously satisfied.  If not, then writing $r = \frac{\kappa_2}{\kappa_1} \geq 1$, the conditions become
 $$r \leq 1 + \frac{2 \Phi'}{\Phi'' F} \mbox{ and } r \geq \frac{1}{1+ \frac{2 \Phi'}{\Phi'' F}} \mbox{.}$$
The second if these is obviously true since $\Phi', \Phi'', F > 0$, while the first is a curvature pinching condition.  Suppose there is a first time $t \geq 0$ where equality is achieved in the first condition above ($t$ could be $0$, but if the inequality is strict on $M_0$ then we will have $t>0$).  At this time, the maximum of $G$ is nonincreasing and thus the pinching ratio does not deteriorate.  We conclude that
  $$r \leq 1 + \frac{2 \Phi'}{\Phi'' F}$$
  continues to hold as long as the solution to \eqref{E:theflow} exists.  That the pinching ratio is strictly decreasing unless $M_t$ is a sphere follows by the strong maximum principle: if there were a time $t_0>0$ where a new maximum of $G$ were attained, then $G$ would be identically constant.  But any surface has an umbilic point, so there is a point where $G=0$ and thus $G\equiv 0$.  It follows that $M_t$ is a sphere.
  
In the axially symmetric case, $\nabla_2 h_{11} \equiv 0$ (see, for example, \cite[Lemma 3.2]{MMW}) but there is no symmetry in $\kappa_1$ and $\kappa_2$.  We require that at the maximum of $G$,  
  $$2 \Phi' \kappa_2 + \Phi'' F \left( \kappa_2 - \kappa_1 \right) \geq 0 \mbox{.}$$
  If it happens that $\Phi'' = 0$ at the maximum of $G$, then the above is clearly satisfied.  Otherwise, if $\kappa_2 \geq \kappa_1$ at the maximum of $G$, then the above is always true, reading
$$\frac{\kappa_2}{\kappa_1} \geq \frac{1}{1+ \frac{2 \Phi'}{\Phi'' F}} \mbox{.}$$
If on the other hand $\kappa_1 > \kappa_2$ at the maximum of $G$, then the above reads
$$\frac{\kappa_1}{\kappa_2} \leq \frac{1}{1+ \frac{2 \Phi'}{\Phi'' F}} \mbox{.}$$
The same argument as in the $n=2$ case gives that with this initial requirement on the pinching ratio, the maximum of $G$ does not deteriorate and thus the pinching ratio does not deteriorate.  To see that the pinching ratio is strictly decreasing, we again appeal to the strong maximum principle.  If a new maximum of $G$ is attained at some $\left( x_0, t_0\right)$, $t_0>0$, then $G$ is identically equal to a positive constant.  The evolution equation for $G$ then yields
$$0 \equiv -\frac{2n\left( n-1\right) F}{\kappa_2 H^3} \left[ \Phi'' F \left( \kappa_2 - \kappa_1\right) + 2 \Phi' \kappa_2\right] \left( \nabla_1 h_{22} \right)^2 - \left( \Phi'\, F - \Phi\right) \frac{2n\kappa_1 \kappa_2 \left( \kappa_1 - \kappa_2\right)^2}{H^3} \mbox{.}$$
This is the sum of two nonpositive terms, so each must be identically equal to zero.  If the pinching requirement is actually strict, then the first term and axial symmetry imply $\left| \nabla A \right| \equiv 0$ and thus $M_t$ is a sphere.  Otherwise, from the second term we see that we must have
$$0 \equiv \left( \Phi'\, F - \Phi\right) \frac{2n\kappa_1 \kappa_2 \left( \kappa_1 - \kappa_2\right)^2}{H^3} = \left( \Phi'\, F - \Phi\right) \frac{2 \kappa_1 \kappa_2}{H} C \mbox{.}$$
Since $C>0$ and $M_t$ is strictly convex, it must be that $M_t$ has
$$\Phi'\left( F\left( \mathcal{W} \right) \right)  F\left( \mathcal{W}\right) - \Phi\left( F\left( \mathcal{W} \right) \right) \equiv 0 \mbox{.}$$
If $M_t$ also has $F\left( \mathcal{W}\right)$ identically constant, then by \cite{EH89}, $M_t$ is a sphere.  If not, then the solution of the above ordinary differential equation is $\Phi\left( s\right) = \Phi\left( 1\right) s$.  This particular case of $\Phi$ was covered in \cite{MMW}.  We conclude for the $\Phi$ of interest here that it is not possible for $G$ to be constant unless $M_t$ is a sphere, so the pinching ratio is strictly decreasing unless $M_t$ is a sphere.\hspace*{\fill}$\Box$
\mbox{}\\[8pt]
\noindent \textbf{Remark:} The above pinching requirement is easily seen to be equivalent to that in \cites{Asurf, MMW} for the case $\Phi\left( F \right) = F^\alpha$, $\alpha>1$.  For more general $\Phi$ here it becomes a condition on the initial surface $M_0$.\\[8pt]

\subsection{Case of Conditions \ref{T:Fconds} e) iii) a) and b)}

As in \cite{A94}, Condition \ref{T:Fconds} e) iii) a) implies b), so only b) needs to be considered.  We use the evolution of the same quantity as in \cite{A94} but in this case a condition on $\Phi$ is needed.  As in Lemma \ref{T:evlneqns}, h), we have
\begin{multline} \label{E:HonF}
  \frac{\partial}{\partial t} \left( \frac{H}{F} \right) = \Delta_{\dot \Phi} \left( \frac{H}{F} \right) + \frac{2}{F} \dot \Phi^{kl} \nabla_k F \nabla_l \left( \frac{H}{F} \right) + \frac{\Phi'}{F} \ddot F^{kl, rs} \nabla^i h_{kl} \nabla_i h_{rs}\\
+ \frac{\Phi''}{F^2} \left( F g^{ij} - H \dot F^{ij} \right) \nabla_i F \nabla_j F + \frac{1}{F^2} \left( \Phi' F - \Phi \right) \left( H \dot F^{kl} h_k^{\ p} h_{pl} - F\left| A\right|^2 \right) \mbox{.}
\end{multline}
We wish to show the maximum of $\frac{H}{F}$ is nonincreasing.  The zero order term in \eqref{E:HonF} is nonpositive in view of Conditions \ref{T:Phiconditions}, d) and Lemma \ref{T:Fconvprops}, i).  Using $s\left| \Phi''\right|^2 \leq c\, \Phi'$ we can show that at a maximum of $\frac{H}{F}$ the second and third term, considered together, are nonpositive.  

The second term is a gradient term and the $\ddot F$ term has the right sign for applying the maximum principle since $F$ is concave and $F, \Phi'>0$.  Moreover, from Corollary \ref{T:consign} there is an absolute constant $\overline{C}>0$ such that
$$\frac{\Phi'}{F} \ddot F^{kl, rs} \nabla^{i} h_{kl} \nabla_i h_{rs} \leq -\overline{C} \frac{\Phi'}{F}\frac{\left| \nabla A\right|^2}{\left| A \right|} \mbox{.}$$

We  estimate, using Conditions \ref{T:Phiconditions}, h),
\begin{multline*}
  \frac{\Phi''}{F}\left[ g^{ij} - \left( \frac{H}{F}\right) \dot F^{ij} \right] \dot F^{kl} \dot F^{rs} \nabla_i h_{kl} \nabla_j h_{rs} \leq \frac{\left| \Phi''\right|}{F} \left| g^{ij} - \left( \frac{H}{F}\right) \dot F^{ij} \right| \left| \dot F\right|^2 \left| \nabla A\right|^2\\
  \leq \frac{c\, \Phi'}{F^2} \left| g^{ij} - \left( \frac{H}{F}\right) \dot F^{ij} \right| \left| \dot F\right|^2 \left| \nabla A\right|^2 \mbox{.}
  \end{multline*}
Since $F$ is concave, we always have $\frac{H}{F} \geq n$ by Lemma \ref{T:Fprops}, i).  At an umbilic point $\left( k, \ldots, k\right)$ we have $\frac{H}{F} = n$ and $g^{ij} - \left( \frac{H}{F}\right) \dot F^{ij} = 0$, the $n\times n$ zero matrix.  Let $\tilde c >n$ be the minimal constant such that
\begin{equation} \label{E:coeff}
  n \leq \frac{H}{F} \leq \tilde c \implies \frac{\left| A\right|}{F} \left| g^{ij} - \left( \frac{H}{F}\right) \dot F^{ij} \right| \left| \dot F\right|^2 \leq \frac{\overline{C}}{c} \mbox{.}
  \end{equation}
Then the second and third terms of \eqref{E:HonF} considered together are nonpositive as required.  The weak maximum principle gives that the maximum of $\frac{H}{F}$ is nonincreasing.

Using the strong maximum principle we can show that in fact the maximum of $\frac{H}{F}$ is strictly decreasing unless $M_t$ is a sphere.  If a new maximum of $\frac{H}{F}$ were attained at some $\left( x_0, t_0\right)$, $t_0>0$, then $\frac{H}{F}$ would be identically constant and we would again have 
\begin{multline} \label{E:strong}
  0 \equiv \frac{\Phi'}{F} \ddot F^{kl, rs} \nabla^i h_{kl} \nabla_i h_{rs} + \frac{\Phi''}{F}\left[ g^{ij} - \left( \frac{H}{F}\right) \dot F^{ij} \right] \nabla_i F \nabla_j F \\
  + \frac{1}{F^2} \left( \Phi' F - \Phi \right) \left( H \dot F^{kl} h_k^{\ p} h_{pl} - F\left| A\right|^2 \right) \mbox{.}
  \end{multline}
This is a sum of three nonpositive terms, so each must be identically equal to zero.  Since the concavity of $F$ is strict in nonradial directions, we see from Corollary \ref{T:consign} that $\left| \nabla A \right| \equiv 0$ and thus $M_t$ is a sphere.  We again conclude the pinching ratio is strictly decreasing unless $M_t$ is a sphere.\hspace*{\fill}$\Box$

\subsection{Case of Conditions \ref{T:Fconds} e) ii)}

Here we will use an argument similar to Schulze \cite{Schulze}; the degree-zero quantity $\frac{K}{F^n}$ has been used elsewhere for the case where $F$ is convex (egs \cites{C, A94, M05}).  Under the flow \eqref{E:theflow}, $\frac{K}{F^n}$ evolves according to Lemma \ref{T:evlneqns}, i).  We wish to show the minimum of $\frac{K}{F^n}$ does not decrease under \eqref{E:theflow}.  Using Lemma \ref{T:Fprops}, i) and the fact that $\Phi' F - \Phi \geq 0$, we observe that the zero order term in the evolution equation is nonnegative.  We will also neglect the norm-like term which is nonnegative.  For the third term we use that since $M_t$ is convex, the principal radii of curvature satisfy $r_i= \frac{1}{\kappa_i} \leq \frac{1}{H}$ and thus $\left( h^{-1}\right)^{ij} \geq \frac{1}{H} g^{ij}$, thus we can estimate
$$\Phi' \frac{K}{F^n} \left( h^{-1}\right)^{ij} \ddot F^{kl,rs} \nabla_i h_{kl} \nabla_j h_{rs} \geq \Phi' \frac{K}{F^n} \frac{1}{H}  \ddot F^{kl,rs} \nabla^i h_{kl} \nabla_i h_{rs} \geq \underline{C} \frac{\Phi'}{n} \frac{K}{F^{n+1}}  \frac{\left| \nabla A\right|^2}{\left| A\right|} \mbox{,}$$
where the last inequality holds in view of Lemmas \ref{T:Fprops} i) and \ref{T:consign}.  For the fourth $\Phi''$ term we estimate
\begin{multline*}
  \Phi'' \frac{K}{F^{n+1}} \left[ F\left( h^{-1} \right)^{ij} - n \dot F^{ij} \right] \nabla_i F \nabla_j F \leq \left| \Phi''\right| \frac{K}{F^{n+1}} \left| F\left( h^{-1} \right)^{ij} - n \dot F^{ij} \right| \left| \dot F\right|^2 \left| \nabla A \right|^2 \\
  \leq \frac{c \Phi'}{F} \frac{K}{F^{n+1}} \left| F\left( h^{-1} \right)^{ij} - n \dot F^{ij} \right| \left| \dot F\right|^2 \left| \nabla A \right|^2 \mbox{,}
  \end{multline*}
where we have used Conditions \ref{T:Phiconditions} h).  Therefore we will have
$$\Phi'' \frac{K}{F^{n+1}} \left[ F\left( h^{-1} \right)^{ij} - n \dot F^{ij} \right] \nabla_i F \nabla_j F \leq \underline{C} \frac{\Phi'}{n} \frac{K}{F^{n+1}}  \frac{\left| \nabla A\right|^2}{\left| A\right|}$$
provided the following degree-zero homogeneous inequality holds:
$$\frac{\left| A\right|}{F} \left| F\left( h^{-1} \right)^{ij} - n \dot F^{ij} \right| \left| \dot F \right|^2 \leq \frac{\underline{C}}{n\, c} \mbox{.}$$
Now, the matrix in the above norm is equal to zero on the sphere (on which $\frac{K}{F^n} \equiv 1$) so by continuity there will be a minimal constant $0<\hat C<1$ such that
$$\frac{K}{F^n} \geq \hat C \implies \frac{\left| A\right|}{F} \left| F\left( h^{-1} \right)^{ij} - n \dot F^{ij} \right| \left| \dot F \right|^2 \leq \frac{\underline{C}}{n\, c} \mbox{.}$$
Thus while the curvatures of $M_t$ satisfy $\frac{K}{F^n} \geq \hat C$, we have for almost every $t$
$$\frac{d}{dt} \min_{M_t} \frac{K}{F^n} \geq 0$$
and thus the corresponding curvature pinching is preserved.

To see that in fact the quantity is strictly increasing unless $M_t$ is a sphere we use the strong maximum principle: if $\frac{K}{M_t}$ attained a new minimum at some $\left( x_0, t_0\right)$, $t_0>0$ then $M_t$ would have $\frac{K}{F^n}$ identically constant, and its evolution equation would reduce to the sum of nonnegative terms.  In particular, because of the strict convexity of $F$, $M_t$ would have to have $\left| \nabla A \right| \equiv 0$ and thus $M_t$ would be a sphere.\hspace*{\fill}$\Box$

\subsection{Case of Conditions \ref{T:Fconds} e) v) b)} As in \cite{AM12}, control on the quantity $Z_\sigma = \left| A^0\right|^2 - \sigma H^2$ provides curvature pinching, via Lemma \ref{T:tight}, i).  We use the evolution equation for $Z_\sigma$ similarly as in \cite{AM12}, Section 5.  As usual, some differences arise because we have $\Phi\left( F\right)$ in place of $F^\alpha$.  

Suppose $Z_\sigma \leq 0$ everywhere on $M_0$ and there is a first time $t_0\geq 0$ where $Z_\sigma$ at $\left( x_0, t_0\right)$ is equal to zero.  Of the zero-order terms in the evolution equation for $Z_\sigma$, Lemma \ref{T:evlneqns}, i), the first is obviously equal to zero while the second is nonpositive in view of Condition \ref{T:Phiconditions}, d) (first option) and Lemma \ref{T:tight}, ii).

For the $\ddot \Phi$ term we estimate using Lemma \ref{T:AMHa} and Condition \ref{T:Phiconditions}, h)
\begin{align*}
 & 2 \left[ h^{ij} - \left( \frac{1}{n} + \sigma \right) g^{ij} H \right] \ddot \Phi^{kl, rs} \nabla_i h_{kl} \nabla_j h_{rs}\\
 &= 2\left( \Phi' \ddot F^{kl, rs} + \Phi'' \dot F^{kl} \dot F^{rs} \right) \left[ h^{ij} - \left( \frac{1}{n} + \sigma \right) g^{ij} \right]\nabla_i h_{kl} \nabla_j h_{rs}\\
 &\leq 2 \left( \Phi' \left| \ddot F\right| + \left| \Phi'' \right| \left| \dot F\right|^2 \right)\left| h^{ij} - \left( \frac{1}{n} + \sigma \right) g^{ij} \right| \left| \nabla A\right|^2\\
 &\leq 2 \left[ \Phi' \mu + \left| \Phi''\right| H \left( 1 + \mu \frac{\left| A^0\right|}{H} \right)^2 \right] \sqrt{\sigma\left( 1 + n \sigma \right)}  \left| \nabla A\right|^2\\
 &\leq 2 \left[ \Phi' \mu + F \left| \Phi'' \right| \frac{2}{2 - \mu \sigma} \left( 1 + \mu \sqrt{\sigma} \right)^2\right] \sqrt{\sigma \left( 1 + n \sigma \right)} \left| \nabla A\right|^2\\
 &\leq 2 \, \Phi' \left[ \mu +  \frac{2\, c}{2 - \mu \sigma} \left( 1 + \mu \sqrt{\sigma} \right)^2\right] \sqrt{\sigma \left( 1 + n \sigma \right)} \left| \nabla A\right|^2 \mbox{.}
\end{align*}
  
Next we estimate
\begin{align*}
 & \dot \Phi^{ij} \left( \nabla_i h_{kl} \nabla_j h^{kl} - \frac{1}{n} \nabla_i H \nabla_j H \right)\\
 & = \Phi' \dot F^{ij} \left( \nabla_i h_{kl} - \frac{1}{n} g_{kl} \nabla_i H \right) \left( \nabla_j h^{kl} - \frac{1}{n} g^{kl} \nabla_j H \right)\\
  &\geq \Phi' \left( 1 - \mu \frac{\left| A^0\right|}{H} \right) \left| \nabla A^0 \right|^2 \geq \frac{2\left( n-1\right)}{3n} \left( 1- \mu\sqrt{\sigma} \right) \Phi' \left| \nabla A \right|^2 
 \end{align*}
 and, using Lemma \ref{T:gradA},
 \begin{multline*}
2\sigma \dot \Phi^{ij} \nabla_i H \nabla_j H = 2\sigma \Phi' \dot F^{ij} \nabla_i H \nabla_j H\\
 \leq 2 \sigma \left( 1 + \mu \sqrt{\sigma} \right) \Phi' \left| \nabla H \right|^2 
\leq \frac{2\left( n+2\right)}{3} \sigma \left( 1 + \mu \sqrt{\sigma} \right) \Phi' \left| \nabla A \right| \mbox{.}
\end{multline*}

Putting these estimates together we see that the gradient terms of the evolution equation for $Z_\sigma$ may be estimated above by
\begin{multline} \label{E:gt}
  \Phi' \left\{ - \frac{4\left( n-1\right)}{3n} \left( 1- \mu \sqrt{\sigma} \right) + 2\sqrt{\sigma\left( 1 + n\sigma\right)} \left[ \mu + \frac{2\, c}{2 - \mu \sigma} \left( 1 + \mu\sqrt{\sigma}\right)^2 \right] \right.\\
  \left. + \frac{2\left( n+2\right)}{3} \sigma \left( 1 + \mu \sqrt{\sigma} \right)\right\} \left| \nabla A\right|^2 \mbox{.}
  \end{multline}

Recall that $\sqrt{\sigma} < \frac{1}{\mu}$ ensures $\dot F$ is positive definite.  Taking $\sqrt{\sigma}$ smaller still, the first term in the above bracket is negative and enough to ensure overall the expression is nonpositive for $0\leq \sigma \leq \sigma_0$ for some $\sigma_0>0$ depending on $n$ and $\mu$.

We conclude that provided $0\leq \sigma \leq \delta_0$, $Z_\sigma$ does not increase, and thus via Lemma \ref{T:tight}, i) curvature pinching is maintained under the flow \eqref{E:theflow}.\\

To see that curvature pinching is strictly improving in this case we also follow \cite{AM12}.  Set $\overline{h} = \sup_{M_0} H$ and consider for $\lambda>0$ to be chosen the evolution of $Z_\lambda = \left| A^0\right|^2 - \sigma_0 \overline{h}^\lambda H^{2-\lambda}$.  We show that $Z$ cannot attain a new zero maximum at a point with $H\geq \overline{h}$, concluding that, under \eqref{E:theflow},
\begin{equation} \label{E:A0improve}
  \left| A^0\right|^2 \leq \min \left\{ \sigma_0 H^2, \sigma_0 \overline{h}^\lambda H^{2-\lambda} \right\} \mbox{.}
  \end{equation}
Setting $\sigma = \sigma_0 \left( \frac{\overline{h}}{H} \right)^\lambda \leq \sigma_0$, using Lemma \ref{T:evlneqns}, the evolution equation for $Z_\lambda$ is
\begin{align} \label{E:Zl}
   \frac{\partial}{\partial t}  Z_\lambda 
   &= \Delta_{\dot \Phi} Z_\lambda + 2\left[ h^{ij} - \left( \frac{1}{n} + \sigma \right) H g^{ij} \right] \ddot \Phi^{kl, rs} \nabla_i h_{kl} \nabla_j h_{rs}\\ \nonumber
&\quad -2\dot \Phi^{ij} \left( \nabla_i h_{kl} \nabla_j h^{kl} - \frac{1}{n} \nabla_i H \nabla_j H \right)+ 2\sigma \dot \Phi^{ij} \nabla_i H \nabla_j H + \lambda\, \sigma\, H \ddot \Phi^{kl, rs} \nabla^i h_{kl} \nabla_i h_{rs}\\ \nonumber
& \quad - \lambda \left( 3-\lambda \right) \sigma \dot \Phi^{ij} \nabla_i H \nabla_j H + 2  \dot \Phi^{kl} h_k^{\ p} h_{pl} Z_\lambda + \lambda \sigma \dot \Phi^{kl} h_{k}^{\ m} h_{ml} H^2\\ \nonumber
  &\quad - \frac{2}{n} \left( \Phi' F - \Phi \right) \left[ n \, C - \left( 1 + n\sigma\right) H \left| A\right|^2 + \frac{n\, \lambda\, \sigma}{2} H \left| A \right|^2 \right] \mbox{.} \nonumber
\end{align}
Choosing now $\lambda \leq 3$ we can discard the above $\left( 3-\lambda \right)$ term.  Using Condition \ref{T:Phiconditions} d), the other zero order terms are bounded by
\begin{multline*}
  \lambda \, \sigma \Phi' F H^3 - \frac{2}{n} \left( \Phi' F -\Phi \right) \sigma \left( 1 + n \sigma \right) \left( 1 - \sqrt{n\left( n-1\right) \sigma} \right) H^3\\
   \leq \Phi' \left[ \lambda \, \sigma - \frac{2}{n} \varepsilon \sigma  \left( 1 + n \sigma \right) \left( 1 - \sqrt{n\left( n-1\right) \sigma} \right) \right] F H^3 \mbox{,}
   \end{multline*}
which is nonpositive for sufficiently small $\lambda$.

The gradient terms may be estimated as for \eqref{E:gt}; with the extra gradient term in \eqref{E:Zl} we now have
\begin{multline*} 
  \Phi' \left\{ - \frac{4\left( n-1\right)}{3n} \left( 1- \mu \sqrt{\sigma} \right) +\left( 2\sqrt{\sigma\left( 1 + n\sigma\right)} + \lambda \sigma \right) \left[ \mu + \frac{2\, c}{2 - \mu \sigma} \left( 1 + \mu\sqrt{\sigma}\right)^2 \right] \right.\\
  \left. + \frac{2\left( n+2\right)}{3} \sigma \left( 1 + \mu \sqrt{\sigma} \right)\right\} \left| \nabla A\right|^2 \mbox{.}
  \end{multline*}
Taking $\lambda$ smaller again if needed, we can ensure this term is also negative.  We conclude the corresponding $Z_\lambda$ cannot attain a new zero maximum and thus \eqref{E:A0improve} holds.\hspace*{\fill}$\Box$
  
  \subsection{Case of Conditions \ref{T:Fconds} e) iv)}
  Here we will again use Andrews' generalisation of Hamilton's tensor maximum principle, but this time with the parametrisation of the evolving hypersurface via its support function.  As a generalisation of the tensor used in \cite[Lemma 11]{AMZ11}, where the speed was degree-one homogeneous, here we instead set $T_{ij}=r_{ij} - C \Phi^{-1} \overline{g}_{ij}$, where we note $\Phi$ is a function of $F_*^{-1}$ and we choose $C>0$ such that initially $T_{ij} \leq 0$.  Using the evolution equations from Lemma \ref{T:evlneqns2} we have
\begin{multline} \label{E:Tij}
  \frac{\partial}{\partial t} T_{ij} = \Phi' F_*^{-2} \dot F_*^{pq} \overline{\nabla}_p \overline{\nabla}_q T_{ij} + \Phi' F_*^{-2} \ddot F_*^{pq, rs} \overline{\nabla}_i r_{pq} \overline{\nabla}_j r_{rs} - \left( 2\Phi' + \frac{\Phi''}{F_*} \right) F_*^{-3} \overline{\nabla}_i F_* \overline{\nabla}_j F_*\\
  + \frac{2 C \Phi'}{\Phi^3} F_*^{-2} \dot F^{pq}_* \overline{\nabla}_p \Phi \overline{\nabla}_q \Phi \overline{g}_{ij} - \Phi' F_*^{-2} \mbox{\textnormal{tr} } \dot F_* T_{ij} + \left( \Phi' F_*^{-1} - \Phi \right) \overline{g}_{ij} \mbox{.}
  \end{multline}
  Suppose there is a first $\left( x_0, t_0\right)$ where $T_{ij}$ has a null eigenvector, that is, there exists a nonzero vector $\xi$ such that $T_{ij}\left( x_0, t_o\right) \xi^i=0$; equivalently $r_{ij}\xi^i = C\Phi^{-1} \xi_j$.  We may assume all eigenvalues are distinct.  Let us choose coordinates such that $\left( r_{ij}\right) = \mbox{\textnormal{diag}}\left( r_1, r_2, \ldots, r_n \right)$ where $r_1$ is the largest principal radii of curvature and $\xi = e_1$, the first coordinate direction.  Applying the null eigenvector to the zero order terms of \eqref{E:Tij} at $\left( x_0, t_0\right)$, we have
  $$\left( \Phi' F_*^{-1} - \Phi \right) \left| \xi\right|^2 \leq 0$$
  by the condition on $\Phi$.  For the gradient terms, we need to show that for all totally symmetric three-tensors satisfying $\left( B_{kij} - C \Phi^{-2} \Phi' F_*^{-2} \delta_{ij} \dot F_*^{pq} B_{kpq} \right) \xi^i \xi^j$ we have for some choice of $\Gamma_i^{\ j}$ the inequality
\begin{multline*} 
  \left[ \ddot F_*^{kl, pq} - \frac{1}{F_*}\left( 2 + \frac{\Phi''}{\Phi' F_*} \right) \dot F_*^{kl} \dot F_*^{pq} \right] B_{ikl} B_{jpq} \xi^i \xi^j + \frac{2C\,F_*^{-4}}{\Phi^3} \left( \Phi'\right)^2 \dot F_*^{pq} \dot F_*^{rs} \dot F_*^{ij} B_{ipq} B_{jrs}\left| \xi\right|^2\\
  + 2 \dot F_*^{kl}\left[ 2 \Gamma_{k}^{\ p} \left( B_{lip}- C \delta_{ip} \dot F_*^{rs} B_{lrs}\right) \xi^i - \Gamma_k^{\ p} \Gamma_l^{\ q} T_{pq}\right] \leq 0 \mbox{.}
  \end{multline*}
  In view of our chosen coordinates, this can be rewritten as
  \begin{multline} \label{E:nec}
  \ddot f_*^{kp} B_{1kk} B_{1pp} + 2 \sum_{k<p} \frac{\dot f_*^k - \dot f_*^p}{r_k - r_p} B_{1kp}^2 - \frac{1}{F_*}\left( 2 + \frac{\Phi''}{\Phi' F_*} \right) \dot f_*^k \dot f_*^p B_{1kk} B_{1pp}\\
  + \frac{2C\, F_*^{-4}}{\Phi^3} \left( \Phi'\right)^2 \dot f_*^p \dot f_*^r \dot f_*^i B_{ipp} B_{irr} + 2 \dot f_*^k \left[ 2 \Gamma_k^p \left( B_{k1p} - C\delta_{1p} \dot f_*^r B_{krr}\right) - \Gamma_k^{\ p} \Gamma_k^{\ q} T_{pq} \right] \leq 0 \mbox{.}
  \end{multline}
   
  In our coordinates, the condition on $B$ (coming from the first order derivative condition at the extremum) becomes $B_{k11} = C \Phi^{-2} \Phi' F_*^{-2} \dot f_*^i B_{kii}$ for each $k$.  Moreover, $T_{ij} = t_i\, \delta_{ij}$ (where $t_i = r_i - C\, \Phi^{-1}$), so \eqref{E:nec} becomes
  \begin{multline*}
    \ddot f_*^{kp} B_{1kk} B_{1pp} + 2 \sum_{k>1} \frac{\dot f_*^1 - \dot f_*^k}{r_1 - r_k} B_{k11}^2 + 2 \sum_{1<k<p} \frac{\dot f_*^k - \dot f_*^p}{r_k - r_p} B_{1kp}^2 - \frac{1}{F_*}\left( 2 + \frac{\Phi''}{\Phi' F_*} \right) \frac{B_{111}^2}{C^2 \Phi^{-4} F_*^{-4} \left( \Phi'\right)^2}\\
    + \frac{2C}{\Phi^3} \dot f_*^i \frac{B_{i11}^2}{C^2 \Phi^{-4}} + 2 \dot f_*^k \sum_{p>1} \left[ 2 \Gamma_k^p \left( B_{k1p} - C\delta_{1p}  \frac{B_{k11}}{C\Phi^{-2} F_*^{-2} \Phi'}\right) - \left( \Gamma_k^{\ p} \right)^2 t_{p} \right] \leq 0
  \end{multline*}
  
  Choosing now $\Gamma_k^{\ 1}=0$ for all $k$ (as in \cite{Apinching}), the above may be rewritten as
  \begin{multline*}
    \ddot f_*^{kp} B_{1kk} B_{1pp} + 2 \sum_{k>1} \frac{\dot f_*^1 - \dot f_*^k}{r_1 - r_k} B_{k11}^2 + 2 \sum_{1<k<p} \frac{\dot f_*^k - \dot f_*^p}{r_k - r_p} B_{1kp}^2 - \frac{1}{F_*}\left( 2 + \frac{\Phi''}{\Phi' F_*} \right) \frac{B_{111}^2}{C^2 \Phi^{-4} F_*^{-4} \left( \Phi'\right)^2}\\
    + \frac{2\, \Phi}{C} \dot f_*^i B_{i11}^2 + 2\dot f_*^k  \sum_{p>1} \left[ 2 \Gamma_k^p B_{k1p} - \left( \Gamma_k^{\ p} \right)^2 \left( r_p - r_1 \right) \right] \leq 0
    \end{multline*}
    and as in \cite{Apinching} the last summation terms may be rewritten to give equivalently
      \begin{multline} \label{E:nec2}
    \ddot f_*^{kp} B_{1kk} B_{1pp} + 2 \sum_{k>1} \frac{\dot f_*^1 - \dot f_*^k}{r_1 - r_k} B_{k11}^2 + 2 \sum_{1<k<p} \frac{\dot f_*^k - \dot f_*^p}{r_k - r_p} B_{1kp}^2 - \frac{1}{F_*}\left( 2 + \frac{\Phi''}{\Phi' F_*} \right) \frac{B_{111}^2}{C^2 \Phi^{-4} F_*^{-4} \left( \Phi'\right)^2}\\
    + \frac{2\, \Phi}{C} \dot f_*^i B_{i11}^2 - 2 \sum_{p>1}\frac{\dot f_*^k}{r_1-r_p} B_{1kp}^2 + 2  \sum_{p>1} \dot f_*^k \left( r_1 - r_p \right) \left[ \Gamma_k^p + \frac{B_{1kp}}{r_1 - r_p} \right]^2 \leq 0 \mbox{.}
    \end{multline}
The optimal choice of $\Gamma_k^{\ p}$ for $p>1$ is clearly $\Gamma_k^{\ p} = -  \frac{B_{1kp}}{r_1 - r_p}$.  Turning our attention now to the other terms above, the $\ddot f_*$ term is clearly nonpositive since $f_*$ is concave.  The coefficient of $B_{111}^2$ is
\begin{align*}
&-\frac{1}{F_*}\left( 2 + \frac{\Phi''}{\Phi' F_*} \right) \frac{1}{C^2 \Phi^{-4} F_*^{-4} \left( \Phi'\right)^2} + \frac{2\Phi}{C} \dot f_*^1\\
&= \frac{2\Phi^2}{C^2} \left[ - \left( \frac{\Phi^2 F_*^{2}}{\left( \Phi'\right)^2} + \frac{1}{2} \frac{\Phi'' \Phi^2 F_*}{\left( \Phi'\right)^3} \right) F_* + r_1 \dot f_*^1 \right]\\
&=\frac{2\Phi^2}{C^2} \left[ - \left( \frac{\Phi^2 F_*^{2}}{\left( \Phi'\right)^2} + \frac{1}{2} \frac{\Phi'' \Phi^2 F_*}{\left( \Phi'\right)^3} \right) F_* + F_* - \sum_{k=2}^{n} \dot f_*^k r_k \right] \mbox{.}
\end{align*}
A sufficient condition for the above to be nonpositive is
$$ - \left( \frac{\Phi^2 F_*^{2}}{\left( \Phi'\right)^2} + \frac{1}{2} \frac{\Phi'' \Phi^2 F_*}{\left( \Phi'\right)^3} \right) + 1\leq 0\mbox{,}$$
which is true in view of the second order condition on $\Phi$ in the theorem (Conditions \ref{T:Phiconditions} i).

That the remaining coefficients of $B_{k11}^2$, of $B_{1kk}^2$ and of $B_{1kl}^2$, $1< k <l$ in \eqref{E:nec2} are each nonpositive follows exactly as in the proof of \cite[Lemma 11]{AMZ11}. 

The result follows.  Since the inequality of Conditions \ref{T:Phiconditions} is strict, then the quantity 

$\sup_{v\in T_z \mathbb{S}^n, \left\| v\right\|=1} \left( \left. r\left( v, v\right)\right|_{\left( z, t\right)} \Phi\circ F_*\left( r_{ij}\left( z, t\right) \right)\right)$ is strictly decreasing unless $M_t$ is a totally umbillic sphere, by the same strong maximum principle argument as in \cite{AMZ11}, using also here the strict monotonicity of $\Phi$.

 \section{Long time existence and convergence to a point} \label{S:T}
\newtheorem{Tsoest}{Proposition}[section]

The argument for an upper speed bound under the flow \eqref{E:theflow} goes back to an idea of Chou \cite{Tso}.  The version used here is essentially as in \cite{LL20}.

\begin{Tsoest} \label{T:Tso}
Under the flow \eqref{E:theflow}, while the inradius of $M_t$ is positive, the speed satisfies
$$\Phi\left( F\right) \leq C$$
everywhere on $M_t$.
\end{Tsoest}

\noindent \textbf{Idea of the proof:} Curvature pinching of the previous section implies there exists $\underline{C}>0$ such that
$$\dot F^{kl} h_k^{\ m} h_{ml} \geq \underline{C} F^2 \mbox{.}$$
From Lemma \ref{T:evlneqns}, k) we have, while $u- \delta \geq \frac{ \delta}{2}$ say, for small $\delta >0$,
$$\frac{\partial}{\partial t} Z \leq \Delta_{\dot \Phi} Z + \frac{2}{u-\delta} \left< \nabla u, \nabla Z \right>_{\dot \Phi} + Z^2  + \left( 1 - \delta \underline{C} F \right) \frac{Z\Phi' F}{u-\delta} \mbox{.}$$
We only need to consider the situation where $F> \frac{2}{\delta \underline{C}}$ because if $F\leq \frac{2}{\delta \underline{C}}$ then via monotonicity of $\Phi$ we have 
$$\Phi\left( F\right) \leq \Phi\left( \frac{2}{\delta \underline{C}}\right)$$
which is bounded, and hence $Z$ is bounded.  When $F> \frac{2}{\delta \underline{C}}$ we have
\begin{multline*}
  \frac{\partial}{\partial t} Z \leq \Delta_{\dot \Phi} Z + \frac{2}{u-\delta} \left< \nabla u, \nabla Z \right>_{\dot \Phi} + Z^2  - \frac{1}{2} \delta \underline{C} \frac{Z\Phi' F^2}{u-\delta}\\
   =   \Delta_{\dot \Phi} Z + \frac{2}{u-\delta} \left< \nabla u, \nabla Z \right>_{\dot \Phi} + Z^2 \left( 1 - \frac{\underline{C}\, \delta}{2} \frac{\Phi' F^2}{\Phi} \right)\mbox{.} 
   \end{multline*}
From Conditions \ref{T:Phiconditions}, f) observe that if $F$ is large enough we will have
$$1 - \frac{\underline{C}\, \delta}{2} \frac{\Phi' F^2}{\Phi} \leq -1$$
say in which case $Z$ does not increase.  (Again, if $F$ is not large enough then $\Phi$ is bounded anyway.). The result follows. \hspace*{\fill}$\Box$\\[8pt]

The arguments for regularity of solutions to \eqref{E:theflow} while the inradius is positive are now very similar to those in previous work.  Since $F$ is degree-one homogeneous and $\Phi'\left( F\right)$ has a positive lower bound, evolution equations for $F$ and first derivatives of the graph height function are each uniformly parabolic, so a result of Krylov and Safonov \cite{KS} gives H\"{o}lder continuity of these quantities.  Where $F$ is convex or concave, $C^{2, \alpha}$ regularity follows from a classical result of Krylov \cite{K}, in the cases of $n=2$ or axial symmetry, this regularity follows instead as in \cite{Asurf} and \cite{MMW} respectively, using \cite{APDE}, while in the remaining strong curvature pinching case arguments as in \cite{AM12} apply.  Higher regularity then follows by Schauder estimates and a standard contradiction to the maximal time argument using short time existence gives that $\left| A\right|^2 \rightarrow \infty$ as $T\rightarrow \infty$.  Curvature pinching implies all curvatures are unbounded as $t\rightarrow T$ and convergence to a point follows.\hspace*{\fill}$\Box$

\section{Convergence of the rescaled solution to the unit sphere} 

In order to study the asymptotic shape of $M_t$ towards the final time $T$, as in previous work we will rescale by the time-dependent factor determined by the behaviour of shrinking spheres under the flow.  Like in \cite{CT}, for general $\Phi$ this factor cannot usually be written down explicitly but for a given $\Phi$ and initial hypersurface $M_0$ we will denote by $\Theta\left( t\right)$ the radius of the sphere that shrinks to a point precisely at time $T$, the extinction time for the flow \eqref{E:theflow} with initial hypersurface $M_0$.  We then use the rescaled time parameter $\tau = - \ln \Theta\left( t\right)$ and examine the rescaled solution hypersurfaces $\tilde M_\tau$ given by $\tilde X \left( x, \tau \right) = \frac{1}{\Theta\left( t\right)} \left[ X\left( x, t\right) - p \right]$, where $p$ is the extinction point for the solution to \eqref{E:theflow} with initial data $M_0$.  The rescaled hypersurfaces satisfy the evolution equation
$$\frac{\partial \tilde X}{\partial \tau} = - \Phi \left( \frac{1}{\Theta} \right) \Phi\left( F \right) \tilde \nu + \tilde X \mbox{,}$$
where the outer unit normal $\tilde \nu = \nu$ is unchanged under rescaling.

Our argument proceeds as in \cite{A94} Section 7, without writing down $\Theta\left( t\right)$ explicitly.  As there, curvature pinching implies uniform bounds above and below on the inradius and circumradius of $\tilde M_\tau$ for all $\tau$ and via Theorem \ref{T:Tso} the flow speed of the rescaled hypersurfaces is uniformly bounded above.  In most cases the structure of $\Phi$ implies that $\Phi'$ is uniformly bounded below by a positive constant, so, since $F$ is degree-one homogeneous the equation for the rescaled speed is uniformly parabolic with bounded, measurable coefficients and thus the Harnack inequality gives a positive uniform lower bound on the rescaled speed.\footnote{Of the example speeds given in Section \ref{S:examplespeeds}, those of example (1) have $\Phi'\left( s\right) \geq c_1>0$ in the natural case that $k_1=1$; example (2) has $\Phi'$ bounded below in view of the uniform upper bound on the rescaled speed and example (3) has $\Phi' \geq 1$.  If $\Phi'$ is not uniformly bounded below then additional arguments are needed.  These follow \cite{AM12} in the case of strong curvature pinching or \cite{M17} for the axially symmetric case.}  Regularity of the rescaled solutions now follows by standard arguments similarly as in the unrescaled case, with each of the sets of conditions.  The strict improvement of the pinching ratio in each case implies subsequential convergence to the sphere; stronger convergence may then be deduced via a linearisation argument.\hspace*{\fill}$\Box$


\begin{bibdiv}
\begin{biblist}

\bib{AS}{article}{
author={Alessandroni, R},
author={Sinestrari, C},
title={Evolution of hypersurfaces by powers of the scalar curvature},
journal={Ann. Sc. Norm. Super. Pisa Cl. Sci. (5)},
volume={9},
date={2010},
number={3}, 
pages={541--571},
}
 
 
  \bib{A94}{article}{
  author={Andrews, B H},
  title={Contraction of convex hypersurfaces in Euclidean space},
  journal={Calc. Var. Partial Differential Equations},
  volume={2},
  date={1994},
  number={2},
  pages={151--171},
  }

\bib{AHarnack}{article}{
  author={Andrews, B H},
  title={Harnack inequalities for evolving hypersurfaces},
  journal={Math. Z.},
  volume={217},
  number={2},
  pages={179--197},
  year={1994},
  }
  
  \bib{AGauss}{article}{
    author={Andrews, B H},
    title={Motion of hypersurfaces by Gauss curvature},
    journal={Pac. J. Math.},
    volume={195},
    number={1}, 
    pages={1--34},
    year={2000},
    }
    
    \bib{AAniso}{article}{
      author={Andrews, B H},
      title={Volume-preserving anisotropic mean curvature flow},
      journal={Indiana U. Math. J.},
      volume={50},
      number={2},
      date={2001},
      pages={783--827},
      }

    \bib{Apinching}{article}{
   author={Andrews, B H},
   title={Pinching estimates and motion of hypersurfaces by curvature
   functions},
   journal={J. Reine Angew. Math.},
   volume={608},
   date={2007},
   pages={17--33},
}

\bib{APDE}{article}{
  author={Andrews, B H},
  title={Fully nonlinear parabolic equations in two space variables},
  eprint={http://arxiv.org/abs/math/0402235},
 }
\bib{Asurf}{article}{
author={Andrews, B H},
title={Moving surfaces by non-concave curvature functions},
journal={Calc. Var. Partial Differential Equations},
volume={39},
date={2010},
number={3--4},
pages={649--657},
}

\bib{AHMWWW14}{article}{
  author={Andrews, B H},
  author={Holder, A B},
  author={McCoy, J A},
  author={Wheeler, G E},
  author={Wheeler, V-M},
  author={Williams, G H},
  title={Curvature contraction of convex hypersurfaces by nonsmooth speeds},
  journal={J. reine angew. Math.},
  date={2017},
  pages={169--190},
  number={727},
  }



    \bib{AM12}{article}{
      author={Andrews, B H},
      author={McCoy, J A},
      title={Convex hypersurfaces with pinched principal curvatures and flow of convex
hypersurfaces by high powers of curvature},
      journal={Trans. Am. Math. Soc.},
      volume={364}, 
      date={2012},
      pages={3427--3447},
      }
      
          \bib{AM16}{article}{
      author={Andrews, B H},
      author={McCoy, J A},
      title={Contraction of convex surfaces by nonsmooth functions of curvature},
      journal={Comm. PDE},
      volume={41},
      number={7},
      pages={1089--1107},
      }
      
    \bib{AMZ11}{article}{
      author={Andrews, B H},
      author={McCoy, J A},
      author={Zheng, Y},
      title={Contracting convex hypersurfaces by curvature},
      journal={Calc. Var. Partial Differential Equations},
      volume={47},
      date={2013},
      number={3--4}, 
      pages={611--665},
      } 

\bib{AK}{article}{
  author={Athanassenas, M},
  author={Kandanaarchchi, S},
  title={Convergence of axially symmetric volume preserving mean curvature flow},
  journal={Pacific J. Math.},
  volume={259},
  date={2012},
  number={1},
  pages={41--54},
  }
  
  \bib{B}{article}{
    author={Baker, C},
    title={The mean curvature flow of submanifolds of high codimension},
    journal={PhD Thesis, Australian National University},
    date={2010},
    }

\bib{BP}{article}{
  author={Bertini, M}, 
  author={Pipoli, G},
  title={Volume preserving non homogeneous mean curvature flow in hyperbolic space},
  journal={Differ. Geom. Appl.},
  volume={54}, 
  pages={448--463},
  date={2017},
  }
  
\bib{BS}{article}{
  author={Bertini, M}, 
  author={Sinestrari, C},
  title={Volume-preserving nonhomogeneous mean curvature flow of convex hypersurfaces},
  journal={Ann. Mat. Pura Appl.},
  volume={197}, 
  pages={1295--1309},
  date={2018},
  }


\bib{CRS}{article}{
  author={Cabezas-Rivas, E},
  author={Sinestrari, C},
  title={Volume-preserving flow by powers of the $m$th mean curvature},
  journal={Calc. Var. Partial Differential Equations},
  volume={38},
  date={2010}, 
  number={3--4}, 
  pages={441--469},
  } 
  
  \bib{Ca}{article}{label={Ca},
    author={Caffarelli, L A},
    title={Interior \emph{a priori} estimates for solutions of fully non-linear equations},
    journal={Ann. Math.},
    volume={130},
    date={1989},
    pages={189--213},
    }
    
    \bib{CS}{article}{
      author={Cannon, J R},
      author={Salman, M},
      title={On a class of nonlinear nonclassical parabolic equations},
      journal={Appl. Anal.},
      volume={85},
      date={2006},
      number={1--3}, 
      pages={23--44},
      } 
      
      \bib{CW}{article}{
  author={Chou, K-S},
  author={Wang, X-J},
  title={A logarithmic Gauss curvature flow and the Minkowski problem}, 
  journal={Ann. Inst. H. Poincar\'{e} -- Anal. Non Lin\'{e}aire}, 
  volume={17},
  number={6},
  date={2000},
  pages={733--751},
  }

\bib{C}{article}{label={Ch},
  author={Chow, B},
  title={Deforming convex hypersurfaces by the $n$th root of the Gaussian curvature},
  journal={J. Differential Geometry},
  volume={22},
  pages={117--138},
  date={1985},
  }
  
  \bib{CG}{article}{
    author={Chow, B},
    author={Gulliver, R},
    title={Aleksandrov reflection and geometric evolution equations I: The $n$-sphere and $n$-ball},
    journal={Calc. Var.},
    volume={4},
    date={1996},
    pages={249--264},
    }
    
    \bib{CTasian}{article}{
    author={Chow, B},
    author={Tsai, D H},
    title={Expansion of convex hypersurfaces by non-homogeneous functions of curvature},
    journal={Asian J. Math.},
    volume={1},
    pages={769--784},
    date={1997},
    }
    
    \bib{CT}{article}{
    author={Chow, B},
    author={Tsai, D H},
    title={Nonhomogeneous Gauss curvature flows},
    journal={Indiana Univ. Math. J.},
    volume={47},
    number={3},
    date={1998},
    pages={965--994},
    }

\bib{Co}{article}{label={Co},
  author={Cordes, H},
  title={\"{U}ber die erste Randwertaufgabe bei quasilinearen Differentialgleichungen
zweiter Ordnung in mehr als zwei Variablen}, 
  journal={Math. Ann.},
  volume={131},
  date={1956}, 
  pages={278--312},
  }

\bib{E}{article}{
  author={Espin, T},
   title={A Pinching Estimate for Convex Hypersurfaces Evolving Under a
Nonhomogeneous Variant of Mean Curvature Flow},
   date={2020},
   status={preprint},
   }
   
   \bib{EH89}{article}{
     author={Ecker, K}, 
     author={Huisken, G},
     title={Immersed hypersurfaces with constant Weingarten curvature},
     journal={Math. Ann.},
     volume={283},
     number={2}, 
     pages={329--332},
     date={1989},
     }
   
   \bib{Ge2}{article}{
     author={Gerhardt, C},
     title={Closed Weingarten hypersurfaces in Riemannian manifolds},
     journal={J. Differential Geom.},
     volume={32},
     pages={299--314},
     date={1996},
     }
   
\bib{GK}{article}{
  author={Guilfoyle, B},
  author={Klingenberg, W},
  title={Parabolic classical curvature flows},
  journal={J. Aust. Math. Soc.},
  volume={104},
  date={2018}, 
  number={3},
  pages={338--357},
  } 



\bib{Ha82}{article}{label={Hm},
  author={Hamilton, R S},
  title={Three-manifolds with positive Ricci curvature},
  journal={J. Differential Geom.},
  volume={17},
  number={2},
  pages={255--306},
  date={1982},
  }
  
  \bib{Hartley}{article}{label={Hl},
    author={Hartley, D},
    title={Motion by mixed volume preserving curvature functions near spheres},
    journal={Pacific J. Math.},
    volume={274},
    date={2015},
    number={2},
    pages={437--450},
    }
  
  
  \bib{Hartman}{book}{label={Hn},
    author={Hartman, P},
    title={Ordinary differential equations},
    edition={2nd},
    series={SIAM classics in applied mathematics},
    publisher={Birkh\"{u}aser},
    place={Boston},
    pages={612+xviii},
    }
  
  
\bib{Hu}{article}{
  author={Huisken, G},
  title={Flow by mean curvature of convex surfaces into spheres},
  journal={J. Differential Geom.},
  volume={20},
  date={1984},
  number={1},
  pages={237--266},
  }
  



\bib{HV}{article}{
  author={Huisken, G},
  title={The volume preserving mean curvature flow},
  journal={J. reine angew. Math.},
  volume={382},
  pages={35--48},
  date={1987},
  }
  
  \bib{K}{article}{
    author={Krylov, N V},
    title={Nonlinear elliptic and parabolic equations of second order},
    publisher={D. Reidel},
    date={1978},
    }
  
  \bib{KS}{article}{
    author={Krylov, N V},
    author={Safonov, M V},
    title={A certain property of solutions of parabolic equations with measurable coefficients},
    journal={Izv. Akad. Nauk},
    volume={40},
    date={1980},
    pages={161--175},
    translation={language={English},
    journal={Math USSR Izv.},
    volume={16},
    date={1981},
    pages={151--164},
    },
    }
    
    \bib{LL20}{article}{
      author={Li, G},
      author={Lv, Y},
      title={Contracting convex hypersurfaces in space form by non-homogeneous curvature function},
      journal={J. Geom. Anal.},
      volume={30},
      date={2020},
      pages={417--447},
      }
    
    \bib{LYW}{article}{
      author={Li, G},
      author={Yu, L},
      author={Wu, C},
      title={Curvature flow with a general forcing term in Euclidean spaces},
      journal={J. Math. Anal. Appl.},
      volume={353},
      date={2009},
      pages={508--520},
      }
    
    \bib{L}{article}{label={Li},
      author={Lieberman, G M},
      title={Second order parabolic partial differential equations},
      publisher={World Scientific},
      place={Singapore},
      date={1996},
      }

\bib{LT}{article}{
  author={Lin, Y-C},
  author={Tsai, D-H},
  title={On a simple maximum principle technique applied to equations on the circle},
  journal={J. Differential Equations},
  volume={245},
  date={2008},
  pages={377--391},
  }
  
  \bib{Lu}{book}{label={Lu},
    author={Lunardi, A},
    title={Analytic semigroups and optimal regularity in parabolic problems},
    publisher={Birkh\"{a}user},
    place={Basel},
    date={1995},
    }

\bib{M03}{article}{
  author={McCoy, J A},
  title={The surface area preserving mean curvature flow},
  journal={Asian J. Math.},
  volume={7},
  date={2003},
  number={1},
  pages={7--30},
  }
  
  \bib{M04}{article}{
    author={McCoy, J A},
    title={The mixed volume preserving mean curvature flow},
    journal={Math. Z.},
    volume={246},
    date={2004},
    number={1},
    pages={155--166},
    }
  

  \bib{M05}{article}{
  author={McCoy, J A},
  title={Mixed volume preserving curvature flows},
  journal={Calc. Var. Partial Differential Equations},
  volume={24},
  date={2005},
  pages={131--154},
  }

    \bib{Mselfsim}{article}{
    author={McCoy, J A},
    title={Self-similar solutions of fully nonlinear curvature flows},
    journal={Ann. Scuola Norm. Sup. Pisa Cl. Sci. (5)},
    volume={10},
    date={2011},
    pages={317--333},
    }
    
    \bib{M17}{article}{
    author={McCoy, J A},
    title={More mixed volume preserving curvature flows},
    journal={J. Geom. Anal.},
    volume={27},
    date={2017},
    pages={3140--3165},
    }
    
    \bib{M20}{article}{
    author={McCoy, J A},
    title={Contracting self-similar solutions of nonhomogeneous curvature flows},
    status={preprint},
    }
    
    \bib{MMW}{article}{
      author={McCoy, J A},
      author={Mofarreh, F Y Y},
      author={Wheeler, V-M},
      title={Fully nonlinear curvature flow of axially symmetric hypersurfaces},
      journal={Nonlinear Differ. Equ. Appl.},
      doi={10.1007/s00030-014-0287-9},
      year={2014},
      number={10},
      pages={1--19},
      }
      
      \bib{MW}{article}{
        author={McCoy, J A},
        author={Wheeler, G E},
        title={A classification theorem for Helfrich surfaces},
        journal={Math. Ann.},
        volume={357},
        pages={1485--1508},
        date={2013},
        }
        
        \bib{N}{article}{
          author={Nirenberg, L},
          title={On a generalization of quasi-conformal mappings and its application to elliptic partial differential equations}, 
          book={
          title={Contributions to the theory of partial differential equations},
    series={Annals of Mathematics Studies}, 
          number={33}, 
          publisher={Princeton University Press}, 
          place={Princeton, N. J.},
          date={1954}, 
          pages={95--100},
          },
          }
        
        
        \bib{Schulze}{article}{label={Sc},
          author={Schulze, F},
          title={Convexity estimates for flows by powers of the mean curvature},
          journal={Ann. Sc. Norm. Super. Pisa Cl. Sci. (5)},
          volume={5},
          date={2006},
          number={2}, 
          pages={261--277},
          } 
          
          \bib{Sin}{article}{label={Si},
            author={Sinestrari, C},
            title={Convex hypersurfaces evolving by volume preserving curvature flows},
            journal={Calc. Var.},
            volume={54},
            date={2015},
            number={2},
            pages={261--277},
            }
            
            \bib{Sm97}{article}{
              author={Smoczyk, K},
              title={Harnack inequalities for curvature flows depending on mean curvature},
              journal={N. Y. J. Math.},
              volume={3},
              date={1997},
              pages={103--118},
              }
        
        \bib{Smoczyk}{article}{
   author={Smoczyk, K},
   title={Starshaped hypersurfaces and the mean curvature flow},
   journal={Manuscripta Math.},
   volume={95},
   date={1998},
   number={2},
   pages={225--236},

}


\bib{Tso}{article}{
  author={Tso, K-S},
  title={Deforming a hypersurface by its Gauss-Kronecker curvature},
  journal={Comm. Pure Appl. Math.},
  volume={38},
  year={1985},
  pages={867--882},
  }
  
  \bib{U}{article}{
    author={Urbas, J},
    title={An expansion of convex hypersurfaces},
    journal={J. Differential Geometry},
    volume={33},
    date={1991},
    pages={91--125},
}

 \end{biblist}
 \end{bibdiv}
\end{document}